\documentclass[11pt,twoside]{amsart}

\numberwithin{equation}{section}
\usepackage{latexsym,amssymb,amsmath,amsthm,amsfonts,amscd}
\usepackage{remark}
\usepackage{enumerate}
\usepackage[all]{xy}
\usepackage{graphicx}
\usepackage{tikz}
\usetikzlibrary{matrix,arrows}
\definecolor{VerdeOlivo}{rgb}{0.3,0.5,0.1}
\definecolor{Magenta}{rgb}{.65,0.15,.2}
\definecolor{Gris}{gray}{0.3}
\newtheorem{Theorem}{Theorem}[section] 
\newtheorem{Lemma}[Theorem]{Lemma} 
\newtheorem{Corollary}[Theorem]{Corollary}

\newtheorem{Proposition}[Theorem]{Proposition} 
\newtheorem{Claim}[Theorem]{Claim} 
 
\theoremstyle{definition}

\newremark{Definition}[Theorem]{Definition}
\newremark{Remark}[Theorem]{Remark}
\newremark{Example}[Theorem]{Example}

\begin{document} 

%===================================% 
%===============Titulo================% 
%===================================% 

\title[On the sandpile group of the cone of a graph]{On the sandpile group of the cone of a graph}

%===================================% 
%===============Autores==============% 
%===================================% 

\author{Carlos A.  Alfaro}
\address{
Departamento de
Matem\'aticas\\
Centro de Investigaci\'on y de Estudios Avanzados del
IPN\\
Apartado Postal
14--740 \\
07000 Mexico City, D.F. 
} 
\email[Carlos A.  Alfaro]{alfaromontufar@gmail.com}
\thanks{The first author was partially supported by CONACyT and the second author was partially supported by SNI}

\author{Carlos E. Valencia}

\email[Carlos E.~Valencia\footnote{Corresponding author}]{cvalencia@math.cinvestav.edu.mx, cvalencia75@gmail.com}


%===================================% 
%==============Resumen==============% 
%===================================% 
\begin{abstract} 
In this article we study the sandpile group of the cone of a graph.
After introducing the concept of uniform homomorphism of graphs we prove that 
every surjective uniform homomorphism of graphs induces an injective homomorphism 
between their sandpile groups.
Also, we establish a relationship between the sandpile group of the cone of the 
cartesian product of graphs an the sandpile group of the cone of their factors.
As an application of these results we obtain an explicit description of a set of generators 
of the sandpile group of the cone of the hypercube.
\end{abstract}

\maketitle

%==================================================================% 
%===========================Introducion============================% 
%==================================================================% 

\section{Introduction} 

The {\it sandpile models} were firstly introduced by Bak, Tang and Wiesenfeld in~\cite{bak87} and~\cite{bak88},
and have been studied under several names in statistical physics, theoretical computer science, algebraic graph theory, and combinatorics.

The {\it abelian sandpile model} of a graph was introduced by Dhar in~\cite{dhar90}, which generalizes the sandpile model of a grid given in~\cite{bak87}.
The abelian sandpile model of Dhar~\cite{dhar90} begins with a connected graph $G=(V,E)$ and a distinguished vertex $s\in V$, called the sink.
Dhar~\cite{dhar90} showed that the set of some configurations (a configurations of $G$ is a vector in $\mathbb{N}^{V\setminus s}$), called {\it recurrent configurations}, with the vertex-by-vertex 
sum as a binary operation forms a finite abelian group, called the sandpile group of $G$.
It follows from Kirchhoff's  Matrix-Tree theorem (see e.g. \cite{biggs93}) that the order of the sandpile group of a graph $G$ is the number of spanning trees of $G$.
Mainly, the abelian sandpile group has been studied under the name of sandpile group, denoted by $SP(G,s)$, and critical group, denoted by $K(G)$. 
It has been also studied under other names, such as Jacobian group, Picard group, dollar game, see for instance~\cite{biggs97,biggs99,lorenzini89,lorenzini91}.

The sandpile group has been completely determined for some family of graphs, 
see for instance~\cite{biggs99,threshold,cartesian,levine,lorenzini89,lorenzini91,merris,musiker,regulartrees}.
The sandpile group of the cartesian product has received special interest, for instance the following cartesian products of graphs it has been determined:
$P_4\times C_n$~\cite{p4cn}, $K_3\times C_n$~\cite{k3cn}, $K_m\times P_n$~\cite{kmpn}, 
$C_4\times C_n$~\cite{wang09}, and $K_m\times C_n$~\cite{corralesthesis,wang09p}.
The abstract structure of the sandpile group has been partially described for the
hypercube~\cite{Bai} and the cartesian product of complete graphs~\cite{cartesian}.
In~\cite{dual} it was proved that the sandpile group of a dual graph $G^*$ is isomorphic to the sandpile group of $G$.
Also, in~\cite{berget09} there are established some relations between the sandpile group of a graph $G$ and the sandpile group of its line graph.
In particular, they proved that if $G$ is non bipartite and regular, then $K({\rm \bf line}(G))$ is completely determined as a function of $K(G)$.
Finally, in~\cite{lorenzini08} a relationship between the eigenvalues and eigenvectors of the Laplacian matrix of a graph and their sandpile group is established. 

Given a natural number $n$, the \textit{$n$-cone} of a graph $G$, denoted by $c_n(G)$, is the graph obtained from $G$ when we add a new vertex $s$ to $G$
and $n$ parallel edges between the new vertex $s$ and all the vertices of $G$.
If $n=1$ we simply write $c(G)$ instead of $c_1(G)$.
In this article we study the sandpile group of the cone of a graph.
In particular, we give a partial description of the sandpile group of the cone of the 
cartesian product of graphs as a function of the sandpile group of the cone of their factors.
Also, we introduce the concept of uniform homomorphism of graphs and prove that 
every surjective uniform homomorphism of graphs induces an injective homomorphism 
between their sandpile groups.
As an application of these two results we obtain an explicit description of a set of generators 
of the sandpile group of the cone of the hypercube of dimension $d$.

A \textit{graph} $G$ is a pair $(V,E)$, where $V$ is a finite set and $E$ is a subset of the set of unordered pair of elements of $V$.
The elements of $V$ and $E$ are called \textit{vertices} and \textit{edges}, respectively. 
If $e=\{x,y\}$, then $x$ and $y$ are \textit{incident} to $e$, $x$ and $y$ are the \textit{ends} of $e$ and $x$ and $y$ are {\it adjacents}.
The \textit{multiplicity} between two vertices $u$ and $v$ of a graph, denoted by $m_{u,v}$, is the number of edges with ends $u$ and $v$.
The \textit{degree} of a vertex $x\in G$, denoted by $d_{G}(x)=d(x)$, is the number of incident edges to $x$.

A graph $G'=(V',E')$ is a \textit{subgraph} of the graph $G=(V,E)$, if $V'\subseteq V$ and $E'\subseteq E$. 
An \textit{induced} subgraph $G[V']=(V',E')$ is a subgraph of $G=(V,E)$ such that every edge $e\in E$ that has its ends in $V'$ is in $E'$. 

The article is organized as follows.
In section $2$ the concepts of graph theory that will be needed in the rest of the article are introduced. 
We also give the combinatorial and algebraic definitions of the sandpile group of $G$ with sink $s_G$.

\medskip

In section $3$ we introduce the concept of uniform homomorphism of graphs.
Let $G$ and $H$ be two graphs and $V\subseteq V(H)$.
A $V$-uniform homomorphism between $G$ and $H$, is a mapping $f:V(G)\rightarrow V(H)$ 
such that for all $x\in V$ and $y\in V(H)$
\[
d_{G[\{u\}\cup S_y]}(u)=m_{x,y} \text{ for all } u \in S_x=f^{-1}(x)
\]
and $f: V(G)\setminus f^{-1}(V)\rightarrow V(H)\setminus V$ is the identity isomorphism.
After introducing the concept of a $V$-uniform homomorphism, we prove the main theorem of this section.

\medskip

{\bf Theorem 3.5.} 
{\it If $f: G\rightarrow H$ is a surjective $V$-uniform homomorphism with $f^{-1}(s_H)=\{s_G\}$ and $s_H\notin V\subset V(H)$ such that $V(H)\setminus V$ is a stable set, 
then the induced mapping $\widetilde{f}: SP(H, s_H)\rightarrow SP(G, s_G)$, given by
\[
\widetilde{f}({\bf c})_v=
\begin{cases}
{\bf c}_{f(v)}  & \text{ if } f(v) \in V, \\
{\rm deg}(f)\cdot c_{f(v)}  & \text{ if } f(v) \notin V,
\end{cases}
\] 
is an injective homomorphism of groups.}

\medskip
 
Section $4$ is devoted to the study of the sandpile group of the cone of the cartesian product of graphs.
Let ${\bf a}\in \mathbb{Z}^{V(G)}$ and ${\bf b}\in \mathbb{Z}^{V(H)}$ be configurations of the cones of $G$ and $H$ respectively.
Taking the cartesian product of configurations as
\[
({\bf a}\Box {\bf b})_{(u,v)}={\bf a}_u+{\bf b}_v  \text{ for all } u\in V(G) \text{ and }  v \in V(H),
\] 
then, ${\bf a}\Box {\bf b}$ is a recurrent configuration of the cone of the cartesian product of $G$ and $H$ 
whenever ${\bf a}$ and ${\bf b}$ are recurrent configurations of $G$ and $H$, respectively.
This definition of the cartesian product of configurations leads to the main result of section $4$.

{\bf Theorem 4.4.} 
{\it If ${\bf e}_H$ is the identity of $SP(c(H), s_{c(H)})$, then the mapping 
\[
\widetilde{\pi}_G: SP(c(G), s_{c(G)})\rightarrow SP(c(G \Box H), s_{c(G \Box H)})
\]
given by $\widetilde{\pi}_G ({\bf a})={\bf a}\Box {\bf e}_H$ is an injective homomorphism of groups.}

Finally, in section $5$ we use an explicit description of the sandpile group of a thick graph with three vertices as well as the results obtained in section 3 and section 4, 
to get a concrete description of a set of generators of the sandpile group of the cone of the hypercube of dimension $d$.
More precisely, if $V(Q_d)=\{v_{\bf a} \, | \, {\bf a}\in \{0,1\}^d\}$ is the vertex set of the hypercube of dimension $d$ and
\[
g_{\beta} (r,t)_{v_{\bf a}}=
\begin{cases}
r & \text{ if } \beta\cdot {\bf a} \text{ is even},\\
t & \text{ if } \beta\cdot {\bf a} \text{ is odd},
\end{cases}
\]
for all $\beta\in \{0,1\}^d$.
Then, 
\[
\widetilde{K}_{\beta}=\{ g_{\beta} (r,t)+ (d-|\beta|){\bf 1}\: | \: 0\leq r,t \leq d \text{ and either } r=|\beta| \text{ or } t=|\beta| \}\subset \mathbb{Z}^{V(Q_d)},
\]
is a set of recurrent configurations of $SP(c(Q_d),s_{c(Q_d)})$ which is a subgroup of $SP(c(Q_{d}),s_{c(Q_d)})$ isomorphic to ${\mathbb Z}_{2|\beta|+1}$.
The next theorem gives a description of the sandpile group of the cone of $Q_d$ gluing all the subgroups $\widetilde{K}_{\beta}$.

{\bf Theorem 5.3.}  
{\it Let $k\geq 0$, $d\geq 1$ be natural numbers and let $c_{2k+1}(Q_d)$ be the $2k+1$-cone of the hypercube $Q_d$.
If $s=V(c_{2k+1}(Q_d))\setminus V(Q_d)$, then
\[
SP(c_{2k+1}(Q_d),s)\cong \bigoplus_{i=0}^{d} \mathbb{Z}_{2i+2k+1}^{\binom{d}{i}}.
\]
Furthermore, $SP(c(Q_{d}),s)=\bigoplus_{\beta \in \{0,1\}^d} \widetilde{K}_{\beta}$.}

The introduction of an extra vertex in the cone's construction is fundamental in order to get a better behavior of the sandpile group.
For instance, in 2003, Jacobson, Niedermaier and Reiner~\cite{cartesian} gave a partial description of the sandpile group of the cartesian product of complete graphs.
In the same year, Bai~\cite{Bai} proved that the number of invariant factors of the hypercube $Q_k$ is $2^{k-1}-1$
and gave a formula for the number of occurrences of $\mathbb{Z}_2$ in the elementary divisor form of the sandpile group of $Q_k$.
However, the full structure of the Sylow $2$-subgroup of the sandpile group of the hypercube is still unknown.

%==================================================================% 
%===========================Preliminaries=========================% 
%==================================================================% 

\section{Preliminaries}
Let $G$ be a graph with $V$ as vertex set and $E$ as edge set.
For simplicity, an edge $e=\{x,y\}$ will be denoted by $xy$. 
The sets of two or more edges with the same ends are called \textit{multiple edges}.
A \textit{loop} is an edge incident to a unique vertex.
A \textit{multigraph} is a graph with multiple edges and without loops. 

A \textit{digraph} $G$ is a pair $(V,E)$, where $V$ is a finite set and $E$ is a subset of the set of ordered pair of elements of $V$.
The elements of $V$ and $E$ are called \textit{vertices} and \textit{arcs}, respectively. 
Given an arc $e=(x,y)$, we say that $x$ is the initial vertex of $e$ and $y$ is the terminal vertex of $e$.
The number of arcs with initial vertex $x$ and terminal vertex $y$ will be denoted by $m_{(x,y)}$.
The \textit{out-degree} of a vertex $x$ of a digraph, denoted by $d_{G}^+(x)$, is the number of arcs with initial vertex $x$.
A vertex $x$ is a {\it sink} if its out-degree is zero.
Moreover, a sink $x$ is a {\it global sink} if for every vertex $y\in G$, there exists a directed path from $y$ to $x$.

Given a multigraph $G$ and a vertex $s$ of $G$, let $b(G,s)$ be the digraph with the same vertex set of $G$ and arc set equal to
\[
E(b(G,s))=\left( \bigcup_{xy\in E(G\setminus s)} \{(x,y),(y,x)\}\right) \cup \left( \bigcup_{xs\in E(G)} \{(x,s)\}\right).
\]
Note that, $b(G,s)$ is a digraph with global sink $s$.

Let $G$ be a digraph, $s$ be a global sink of $G$, and $\widetilde{V}$ the set of non-sink vertices of $G$. 

\subsection{The sandpile group}
There exist several ways to define the sandpile group of a digraph.
In this section we will present a combinatorial and an algebraic definition of the sandpile group.

\paragraph{\it Algebraic description.}
One of the simplest ways to define the sandpile group is by using an algebraic description, known as the critical group.
The {\it Laplacian matrix} of $G$, denoted by $L(G)$, is the matrix of $|V|\times |V|$ given by
\[
L(G)_{u,v}=
\begin{cases}
d^+(u)-m_{(u,u)}  & \text{ if } u=v,\\
-m_{(u,v)} & \text{ otherwise}. 
\end{cases}
\]
The {\it reduced Laplacian matrix}, denoted by $L(G,s)$, is the matrix obtained
from $L(G)$ by removing the row and column $s$.

The {\it sandpile group} of $G$ is the cokernel of $L(G,s)$,
\[
SP(G,s)=\mathbb{Z}^{\widetilde{V}}/{\rm Im}\, L(G,s)^t.
\]

\medskip

Another way to define the sandpile group is by using stable and recurrent configurations.
\paragraph{\it Combinatorial description.}
A \textit{configuration} of $(G,s)$ is a vector ${\bf c}\in \mathbb{N}^{\widetilde{V}}$.
A non-sink vertex $v$ is called \textit{stable} if $d^+(v) > {\bf c}_v$, and otherwise is called unstable. 
Moreover, a configuration is called \textit{stable} if every vertex $v$ in $\widetilde{V}$ is stable.
\textit{Toppling} an unstable vertex $u$ in ${\bf c}$ is performed by decreasing ${\bf c}_u$ by the degree $d^+(u)$, 
and adding the multiplicity $m_{(u,v)}$ to each of the vertices $v$ such that $(u,v)\in E(G)$.
Now, let $\Delta_u=d^+(u)-\sum_{uv\in E} m_{(u,v)}{\bf e}_v$, where ${\bf e}_v$  
is the $v$-th canonical vector with a one in the $v$-th coordinate and zeros elsewhere.
Then, $\Delta_u$ is a row of the reduced Laplacian matrix $L(G,s)$ and
toppling $u$ means to subtract $\Delta_u$ from ${\bf c}$. 

By performing a sequence of topplings, we will eventually arrive at a stable configuration, \cite[Lemma 2.4]{holroyd}. 
See~\cite[Example 2.1]{holroyd} for an example of a digraph without global sink and a configuration that does not stabilizes.
Moreover, the stabilization of a unstable configuration is unique, \cite[Theorem 2.1]{meester}.
The stable configuration associated to ${\bf c}$ will be denoted by $s({\bf c})$.
Then, $s({\bf c})={\bf c}-L(G,s)^t \beta$ for some $\beta\in \mathbb{N}^{\widetilde{V}}$.

Now, let $({\bf c}+{\bf d})_u:={\bf c}_u+{\bf d}_u$ for all $u\in \widetilde{V}$ and ${\bf c}\oplus {\bf d}:=s({\bf c}+{\bf d})$. 
A configuration ${\bf c}$ is \textit{recurrent} if it is stable and there exists a non-zero configuration ${\bf r}$ such that $s({\bf c}+{\bf r})={\bf c}$.
The \textit{sandpile group} of $G$, denoted by $SP(G,s)$, is the set of recurrent configurations with $\oplus$ as binary operation. 

\medskip

Given a multigraph $G$ with a distinguished vertex $s$ their sandpile group is defined by $SP(b(G,s),s)$.

\begin{Theorem}\cite[corollary 2.5]{dual} and \cite[Corollary 2.16]{holroyd}. 
Let $G=(V,E)$ be a multigraph (respectively, digraph) with (respectively, global sink) sink $s\in V$, then $SP(G,s)$ is an abelian group.
\end{Theorem}

One of the simplest ways to check when a configuration of a multigraph is recurrent is given by the following result:

\begin{Theorem}[Burning Algorithm]\cite{dhar90} \label{burning}
A configuration ${\bf c} \in \mathbb{N}^{\widetilde{V}}$ is recurrent if and only if there exist an order 
$u_1,u_2,\cdots, u_n$ of the vertices $\widetilde{V}$ such that if ${\bf c}_1={\bf c}+\sum_{i=1}^n \Delta_{u_i}$, and
\[
{\bf c}_{i}={\bf c}_{i-1}-\Delta_{u_{i-1}}\text{ for all } i=2,\ldots,n,
\] 
then $u_{i}$ is an unstable vertex of ${\bf c}_{i}$ for all $i=1,\ldots,n$ and ${\bf c}={\bf c}_n-\Delta_{u_n}$.
\end{Theorem}

There is a generalization of the burning algorithm for digraphs, know as the script algorithm, see~\cite{speer}

For instance, in the next proposition, we shall describe the sandpile group of the multidigraph
$c(\mathcal{K}_2(r,t))$ with $V =\{s, v_1, v_2\}$ as vertex set, $m_{v_1,s}=1$, $m_{v_2,s}=1$, 
$m_{(v_1,v_2)}=r$ , and $m_{(v_2,v_1)}=t$.
If $r=t$ we simply write $c(\mathcal{K}_2(r))$ instead of $c(\mathcal{K}_2(r,t))$.

\begin{Theorem}\cite[theorem 2.34]{alfarothesis}\label{generators}
If $r \in \mathbb{Z}_+$ and $t \in \mathbb{Z}_+$, then 
\[
SP(c(\mathcal{K}_2(r,t)),s) \cong \mathbb{Z}_{r+t+1}.
\]
Moreover, $SP(c(\mathcal{K}_2(r)),s)=\{ (m,l)\: | \: 0\leq m,l \leq d \text{ and } m=r \text{ or } l=r \}$ with
$(r,r)$ as the identity and $(r,0)$ is a generator of $SP(\mathcal{K}_2(r),s)$ with
\[
k(r,0)=
\begin{cases}
(r-j,r) & \text{ if } k=2j\leq 2r,\\
(r,j) & \text{ if } k=2j+1 \leq 2r+1.
\end{cases}
\]
\end{Theorem}

\medskip

It is known that both descriptions are equivalent in the sense that both descriptions define isomorphic groups, \cite[Corollary 2.16]{holroyd}. .
Is not difficult to see that the structure of the sandpile group does not depend on the sink vertex.
However, the set of recurrent configurations of $G$ depends on the sink.
In this article we are not only interested in the abstract structure of the sandpile group, we are also interested
in the set of recurrent configurations and in the description of the subgroups generated by this recurrent configuration.
We are interested in giving a description of the recurrent configurations because they contain a very nice combinatorial structure and some combinatorial information of the graph.
In general it is easier to describe the abstract structure of the sandpile group than to 
give an explicit description of recurrent configurations and their generated subgroups generated.  
For instance, when $G$ is the grid, in~\cite{BorgneIdentidad}, \cite{identity}, and~\cite{DartoisIdentity}
is given a partial characterization of the recurrent configuration that plays the role of the identity. 
The set of recurrent configuration and their generated subgroup has been described only for a few family of graphs.

In the following, every multigraph will be connected and will have a distinguished vertex $s_G\in V(G)$, called \textit{sink}.
Sometimes we will simply write $s$ instead of $s_G$.
The set of non-sink vertices will be denoted by $\widetilde{V}$.

%================================================================================% 
%==============================Graph homomorphisms ===============================% 
%================================================================================% 

\section{Graph homomorphism and the sandpile group} 

In this section we introduce the concepts of uniform homomorphism and weak homomorphism of graphs.
This concepts are similar to the classical concepts of homomorphism and full homomorphism of graphs.
Also we introduce a directed variant of the uniform homomorphism concept, called directed uniform homomorphism.
In the literature there are several concepts that are either equivalent or similar to the concepts 
of uniform homomorphism, weak homomorphism, and directed uniform homomorphism of graphs.
For instance, in~\cite[chapter 5]{godsil} and~\cite[section 5]{directed}  the concept of an equitable partition of a graph was defined.
This concept of equitable partition is equivalent to the concept of directed uniform homomorphism.
In~\cite[section 5]{bicycles}, Berman defined the concept of divisibility of graphs, 
which is closed related to the concept of weak $V$-uniform homomorphism when $V=V(G)$, 
see remark~\ref{bicycles} for a more precise explanation of this equivalence. 

The concept of uniform homomorphism is useful in order to get an insight of the group structure of the sandpile groups of graphs.
For instance, theorem~\ref{uniformhomeo} says that if $f:G\rightarrow H$ is a surjective $(V(H)\setminus s_H)$-uniform homomorphism,
then the induced mapping $\widetilde{f}:SP(H,s_H) \rightarrow SP(G,s_G)$ is an injective homomorphism of groups;
that is, this mapping sends recurrent configurations to recurrent configurations and is compatible with the group structure.
Theorem 5.7 in~\cite{bicycles} shows an equivalent result to the one in theorem~\ref{uniformhomeo}.
Theorem 6.1 in~\cite{directed} shows an equivalent result to the one in theorem~\ref{directeduniformhomeo}.
In~\cite[section 2]{harmonic}  the concept of harmonic morphism was defined (this concept is different form uniform homomorphism) and
a functor between the category of graphs with harmonic morphisms and the category of abelian groups was studied.
In~\cite{treumann} it is explored a functor from the category of graphs with divisibility to the category of abelian groups, see for instance Proposition 19.
Finally, in~\cite{berget09,durgin,lorenzini91} some functorial results on the category of graphs to the category of abelian groups are proved.
For instance in~\cite[Proposition 2]{lorenzini91} and~\cite[Proposition 21]{treumann}  is proved that: 
if $G$ is a connected graph and $G_k$ is the graph obtained by dividing each edge of $G$ in $k$ edges, 
then there exists a surjective function between the sandpile group of $G_k$ and the sandpile group of $G$.
In~\cite{durgin} is introduced the concept of symmetric configuration and quotient graph are discussed,
more precisely Theorem 2.1 proved that the set of symmetric configurations forms a subgroup of the sandpile group.
In~\cite[Theorems 1.3 and 1.5]{berget09} are established homomorphism between the sandpile group 
of the line graph of a graph $G$ and the sandpile group of $G$ and between the sandpile group 
of the line graph of a graph $G$ and the sandpile group of a subdivision of $G$.

\begin{Definition}
Let $G$, $H$ be multigraphs without loops and $V\subseteq V(H)$. 
A \textit{$V$-uniform homomorphism} of $G$ to $H$, denoted by $f: G\rightarrow H$, 
is a mapping $f:V(G)\rightarrow V(H)$ such that for all $x\in V$ and $y\in V(H)$
\[
d_{G[\{u\}\cup S_y]}(u)=m_{x,y} \text{ for all } u \in S_x=f^{-1}(x)
\] 
and $f: V(G)\setminus f^{-1}(V)\rightarrow V(H)\setminus V$ is the identity isomorphism. 
\end{Definition}

If $f:V(G)\rightarrow V(H)$ is a $V$-uniform homomorphism with $V=V(H)$, then we simply say that 
$f$ is a uniform homomorphism.
In the case of directed multigraphs, we define a 
\textit{directed $V$-uniform homomorphism}
as a mapping
$f:V(G)\rightarrow V(H)$ such that $f: V(G)\setminus f^{-1}(V)\rightarrow V(H)\setminus V$ is the identity isomorphism and
for all $x\in V$ and $y\in V(H)$
\[
d^+_{G[\{u\}\cup S_y]}(u)=m_{(x,y)}
\text{ for all } u \in S_x,
\]
where $d^+_{G}(u)$ is the outdegree of the vertex $u$ in the graph $G$, that is, the number of arcs of $G$ with tail $u$.

If $f: G\rightarrow H$ is a $V$-uniform homomorphism, then $S_x$ is a stable set of $G$ for all $x\in V$ because $H$ has no loops.
Moreover, since $G[S_x\cup S_y]$ is a $m_{x,y}$-regular bipartite graph for all $x \neq y\in V$ and $H[V]$ is connected, then $|S_x|=|S_y|$ for all $x,y\in V$.
The degree of a $V$-uniform homomorphism $f: G\rightarrow H$, denoted by ${\rm deg}(f)$, 
is equal to the cardinality of the set $S_x$ for some $x\in V$.
 
\begin{Proposition}\label{grado}
If $f: G\rightarrow H$ is a $V$-uniform homomorphism and $V(H)\setminus V$ is a stable set, then
\[
d_G(u)=
\begin{cases}
d_{H}(f(u)) & \text{ if }  f(u)\in V,\\
{\rm deg}(f)\cdot d_{H}(f(u)) & \text{ if } f(u)\notin V.
\end{cases}
\]
\end{Proposition}
\begin{proof}
If $u\in f^{-1}(V)$, then
\[
d_G(u)=\sum_{y\in V(H)\setminus f(u)} d_{G[\{u\}\cup S_y]} (u)=\sum_{y\in V(H)\setminus f(u)} d_{H[\{f(u)\}\cup \{y\}]} (f(u))=d_{H}(f(u)).
\]
On  the other hand, since $V(H)\setminus V$ is a stable set, then
\[
d_G(u)=\sum_{v\in S_x, x\in V} d_{G[\{v\}\cup \{u\}]} (u)=\sum_{x\in V} {\rm deg}(f)\cdot d_{H[\{x\}\cup \{f(u)\}]} (f(u))={\rm deg}(f)\cdot d_{H}(f(u)).
\]
when $u\notin f^{-1}(V)$.
\end{proof}

The next proposition gives us an alternative description of a uniform homomorphism.

\begin{Proposition}\label{bipartita}
Let $G$ and $H$ be multigraphs without loops.
Then, $f: G\rightarrow H$ is a uniform homomorphism if and only if
\begin{description}
\item[$(i)$] $S_x=f^{-1}(x)$ is an independent set of $G$ for all $x\in V(H)$,
\item[$(ii)$]  $G[S_x\cup S_y]$ is a $m_{x,y}$-regular bipartite graph for all $x \neq y\in V(H)$.
\end{description}
\end{Proposition}

\medskip

Now, we will introduce the classical definitions of a homomorphism and a full homomorphism of graphs in order to compare 
them with the notion of uniform homomorphism.

Let $G$ and $H$ be multigraphs.
A  \textit{homomorphism (respectively, full homomorphism)} is a mapping
\[
f:V(G)\rightarrow V(H)
\]
such that $f(u)f(v) \in E(H)$ if (respectively, and only if) $uv \in E(G)$.

\medskip

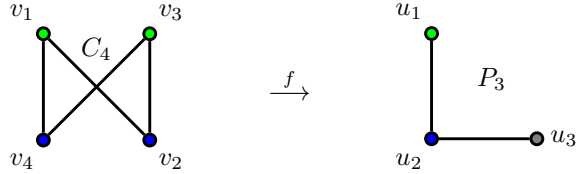
\begin{figure}[h] \centering
	\begin{tabular}{c@{\extracolsep{1cm}}c@{\extracolsep{1cm}}c}
		\begin{tikzpicture}[line width=1.1pt, scale=1]
			\tikzstyle{every node}=[inner sep=0pt, minimum width=4.5pt] 
 			\draw (45:1) node[draw, circle, fill=green] (v1) {};
 			\draw (135:1) node[draw, circle, fill=green] (v3) {};
 			\draw (225:1) node[draw, circle, fill=blue] (v2) {};
 			\draw (315:1) node[draw, circle, fill=blue] (v4) {};
 			\draw (v1) to (v2);
 			\draw (v2) to (v3);
 			\draw (v3) to (v4);
 			\draw (v4) to (v1);
 			\draw (45:1.4) node {\small $v_3$};
 			\draw (135:1.4) node {\small $v_1$};
 			\draw (225:1.4) node {\small $v_4$};
 			\draw (315:1.4) node {\small $v_2$};
			\draw (0,0.5) node {\small $C_4$};
		\end{tikzpicture}
		&
		\begin{tikzpicture}[line width=1.1pt, scale=1]
			\draw (0,0) node {};
			\draw (0,1) node {\small $\overset{f}{\longrightarrow}$};
		\end{tikzpicture}
		&
		\begin{tikzpicture}[line width=1.1pt,scale=1]
			\tikzstyle{every node}=[inner sep=0pt, minimum width=4.5pt] 
 			\draw (0,1.4) node[draw, circle, fill=green] (v1) {};
 			\draw (0,0) node[draw, circle, fill=blue] (v2) {};
 			\draw (1.4,0) node[draw, circle, fill=gray] (v3) {};
 			\draw (v1) to (v2);
			\draw (v2) to (v3);
 			\draw (-.3,1.7) node {\small $u_1$};
 			\draw (-.3,-.3) node {\small $u_2$};
 			\draw (1.76,0) node {\small $u_3$};
			\draw (0.8,0.8) node {\small $P_3$};
		\end{tikzpicture}
	\end{tabular}
	\caption{\small A full homomorphism between $C_4$ and $P_3$.}
	\label{Gfullh}
\end{figure}

The definitions of full homomorphism and isomorphism of graphs are similar. 
The main difference between them is that a full homomorphism is not necessarily bijective; meanwhile an isomorphism is.
By example, let $C_4$ and $P_3$ be graphs as in figure \ref{Gfullh}. 
The mapping $f:V(C_4)\rightarrow V(P_3)$ given by $v_1,v_3 \overset{f}{\mapsto}  u_1$, and $v_2,v_4 \overset{f}{\mapsto}  u_2$ is a full homomorphism.

The following proposition gives us an equivalent way to define a (full) homomorphism of graphs:

\begin{Proposition}\cite[Proposition 1.10 and exercise 10 in page 35]{Nesetril} \label{THEnesetril}
Let $G$ and $H$ be multigraphs without loops.
Then $f: G\rightarrow H$ is an homomorphism if and only if
\begin{description}
\item[$(i)$] $S_x=f^{-1}(x)$ is an independent set of $G$ for all $x\in V(H)$,

\item[$(ii)$] if $xy\notin E(H)$, then $uv\notin E(G)$ for all $u\in S_x$ and $v\in S_y$.
\end{description}
Moreover, $f$ is a full homomorphism if and only if $f$ satisfies conditions $(i)$, $(ii)$, and 
\begin{description}
\item[$(ii')$] if $xy\in E(H)$, then $uv\in E(G)$ for all $u\in S_x$ and $v\in S_y$.
\end{description}
\end{Proposition}

In order to illustrate the concept of uniform homomorphism, 
let $C_3$ and $C_5$ be the cycles with three and five vertices, respectively.
\begin{figure}[h] \centering
		\begin{tikzpicture}[line width=1.1pt, scale=.9]
			\tikzstyle{every node}=[inner sep=0pt, minimum width=4.5pt] 
 			\draw (0,0) {
			+(0:1) node[draw, circle, fill=gray] (v1) {}
 			+(120:1) node[draw, circle, fill=green] (v2) {}
 			+(240:1) node[draw, circle, fill=blue] (v3) {}
 			(v2)+(-1,0) node[draw, circle, fill=green] (v4) {}
 			(v3)+(-1,0) node[draw, circle, fill=blue] (v5) {}
			(v1)+(-4,0) node[draw, circle, fill=gray] (v1p) {}
 			(v1) to (v2) (v4) to (v3) (v1) to (v3) (v5) to (v2) (v4) to (v5)
 			(0:1.4) node {\small $v_1$}
 			(v2)+(0,0.25) node {\small $v_2$}
 			(240:1)+(0,-0.3) node {\small $v_5$}
 			(120:1)+(-1.1,0.3) node {\small $v_4$}
 			(240:1)+(-1.1,-0.3) node {\small $v_3$}
			(v1p)+(-0.3,0) node {\small $v'_1$}
			(-0.0,0) node {\small $C_5$}
			};
			\draw[dashed] (0,0) {
			(v2) to (v3)
			(v1) .. controls +(-90:0.8) and +(-50:1.6).. (v5)
			(v1) .. controls +(90:0.8) and +(50:1.6).. (v4)
			};
			\draw[dashed,gray] (0,0) {
			(v4) to (v1p)
			(v5) to (v1p)
			};
			\draw[dashed,red] (0,0) {
			(v2) to (v1p)
			};
			
			\draw[->, line width=0.5] (2,0) to (3.6,0);
			\draw (2.9,0.3) node {\small $f$};
			
			 \draw (5,0) {
			+(0:1) node[draw, circle, fill=gray] (v1) {}
 			+(120:1) node[draw, circle, fill=green] (v2) {}
 			+(240:1) node[draw, circle, fill=blue] (v3) {}
 			(v1) to (v2) (v2) to (v3) (v1) to (v3)
			(v1)+(0:0.4) node {\small $u_1$}
 			(v2)+(0,0.3) node {\small $u_2$}
 			(v3)+(0,-0.3) node {\small $u_3$}
			+(0.5,0.9) node {\small $C_3$}
			};
			\draw[dashed, gray] (0,0) {
			(v2) to [bend right] (v3)
			};			
		\end{tikzpicture}
	\vspace{-8mm}	
	\caption{\small The mapping $f$.}
	\label{Gfullh1}
\end{figure}
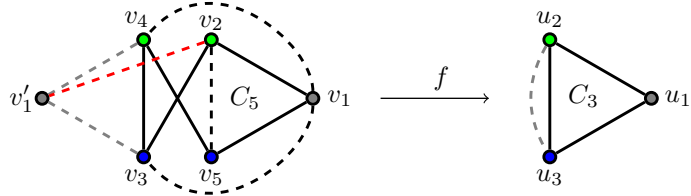
\vspace{-2mm}

The mapping $f:V(C_5)\rightarrow V(C_3)$  given by
\[
\begin{array}{rcc}
v_1 \overset{f}{\longmapsto}  u_1,\\
v_2, v_4 \overset{f}{\longmapsto}  u_2,\\
v_3, v_5 \overset{f}{\longmapsto}  u_3, 
\end{array}
\] 
is a homomorphism of graphs that is neither a full nor uniform homomorphism. 
If we replace $C_5$ by $C_5 + v_2 v_5+v_1v_3+v_1v_4$, we get that the mapping $f$ is a full homomorphism that is not uniform.
Additionally, if we replace $C_5$ by $C_5 + v_2 v_5+v'_1v_4+v'_1v_3$ and $C_3$ by $C_3 + u_2 u_3$,
then the function given by $\widehat{f}(v_i)=f(v_i)$ for all $i=1,\ldots, 5$ and $\widehat{f}(v_1')=u_1$ is a uniform homomorphism,
but $\widehat{f}$ is not a full homomorphism because $v'_1v_2$ is not an edge as required by theorem \ref{THEnesetril}[$(ii)$].

\medskip

The concept of uniform homomorphism of graphs is relevant in the study of the sandpile group of graphs as shown in the following result:

\begin{Theorem}\label{uniformhomeo} 
Let $G$ be a multigraph with sink $s_G$, $H$ be a multigraph with sink $s_H$, $V\subset V(H)$ such that $V(H)\setminus V$ is a stable set and $s_H\notin V$,
$f:G\rightarrow H$ be a surjective $V$-uniform homomorphism such that $f^{-1}(s_H)=\{s_G\}$.
Then the induced mapping $\widetilde{f}:SP(H,s_H) \rightarrow SP(G,s_G)$, given by
\[
\widetilde{f}({\bf c})_v=
\begin{cases}
{\bf c}_{f(v)}  & \text{ if } f(v) \in V, \\
{\rm deg}(f)\cdot {\bf c}_{f(v)}  & \text{ if } f(v) \notin V,
\end{cases}
\] 
is an injective homomorphism of groups, that is, $SP(H, s_H) \triangleleft SP(G, s_G)$.
\end{Theorem}
\begin{proof}
Let $\widehat{f}:\mathbb{Z}^{V(H\setminus s_H)} \rightarrow \mathbb{Z}^{V(G\setminus s_G)}$ be the mapping induced by $f$
given by
\[
\widehat{f}({\bf c})_v=
\begin{cases}
{\bf c}_{f(v)}  & \text{ if } f(v) \in V, \\ 
{\rm deg}(f)\cdot {\bf c}_{f(v)}  & \text{ if } f(v) \notin V.
\end{cases}
\] 
Clearly $\widehat{f}$ is an injective homomorphism of groups.
In order to prove this theorem we need to prove the following facts:
\begin{itemize}
\item If ${\bf c}$ is a recurrent configuration of $(H,s_H)$, then $\widetilde{f}({\bf c})$ is a recurrent configuration of $(G,s_G)$, 
\item $\widetilde{f}({\bf c}_1\oplus {\bf c}_2)=\widetilde{f}({\bf c}_1)\oplus \widetilde{f}({\bf c}_2)$ for all ${\bf c}_1,{\bf c}_2 \in SP(H,s_H)$.
\end{itemize}

Th next claim will be useful to prove this fact.
\begin{Claim}\label{stable}
If ${\bf c}_1$ and ${\bf c}_2$ are configurations of $(H,s_H)$, then 
\[
\widehat{f}(s({\bf c}_1+{\bf c}_2))=s(\widehat{f}({\bf c}_1)+\widehat{f}({\bf c}_2)).\vspace{-1mm}
\]
\end{Claim}
\begin{proof}
By proposition~\ref{grado} a vertex $x\in V(H)\setminus s_H$ can be toppled in the configuration ${\bf c}$ of $(H,s_H)$ 
if and only if the vertices $S_x$ of $G$ can be toppled in the configuration $\widehat{f}({\bf c})$  of $(G,s_G)$. 

On the other hand, since $\widehat{f}(\Delta_x)=\sum_{v\in S_x} \Delta_{v}$ for all $x\in V(H)\setminus \{s_H\}$
and $s({\bf c})={\bf c}-\sum_{w\in W} \Delta_w$ for some multiset $W$ of $V(H)\setminus s_H$,
then 
\begin{eqnarray*}
\widehat{f}(s({\bf c}_1+{\bf c}_2))&=&\widehat{f}\left( {\bf c}_1+{\bf c}_2-\sum_{w\in W} \Delta_w\right)=
\widehat{f}({\bf c}_1)+\widehat{f}({\bf c}_2)-\sum_{w\in W} \widehat{f}(\Delta_w)\\
&=& \widehat{f}({\bf c}_1)+\widehat{f}({\bf c}_2)-\sum_{w\in W} \sum_{v\in S_w} \Delta_{v}=s(\widehat{f}({\bf c}_1)+\widehat{f}({\bf c}_2)). \vspace{-7mm}
\end{eqnarray*}
\end{proof}

Clearly, ${\bf c}$ is a stable configuration of $(H,s_H)$ if and only if $\widehat{f}({\bf c})$ is a stable configuration of $(G,s_G)$.
Furthermore, if ${\bf c}$ is a recurrent configuration of $(H,s_H)$, then there exists a configuration ${\bf u}$ of $(H,s_H)$ such that $s({\bf c}+{\bf u})={\bf c}$.
Thus, by claim~\ref{stable} $s(\widehat{f}({\bf c})+\widehat{f}({\bf u}))=\widehat{f}(s({\bf c}+{\bf u}))=\widehat{f}({\bf c}$) 
and therefore $\widehat{f}({\bf c})$ is a recurrent configuration of $(G,s_G)$.
Finally, $\widetilde{f}({\bf c}_1\oplus {\bf c}_2)=\widetilde{f}(s({\bf c}_1+{\bf c}_2))=s(\widetilde{f}({\bf c}_1)+\widetilde{f}({\bf c}_2))=
\widetilde{f}({\bf c}_1)\oplus \widetilde{f}({\bf c}_2)$
for all ${\bf c}_1,{\bf c}_2 \in SP(H,s_H)$.
\end{proof}

\begin{Remark}\label{contraction}
Note that, a mapping $f:G\rightarrow H$ is a surjective uniform homomorphism if and only if
the induced mapping $\check{f}:\check{G}\rightarrow H$ is a surjective $(V(H)\setminus s_H)$-uniform  homomorphism,
where $\check{G}=G/ f^{-1}(s_H)$ is the graph obtained from $G$ when we contract all the vertices in $f^{-1}(s_H)$ to a single vertex $s_G$.
For instance, consider the next graphs with $f:G\rightarrow H$ given by
\[
\begin{array}{lcr}
u_1, u'_1 & \overset{f}{\longmapsto} & s_H,\\
u_2, u_4 & \overset{f}{\longmapsto}  & v_2,\\
u_3, u_5 & \overset{f}{\longmapsto}  & v_3.
\end{array}
\]
\vspace{-4mm}
\begin{figure}[h]\centering
\begin{tikzpicture}[line width=1.1pt,scale=.9]
			\tikzstyle{every node}=[inner sep=0pt, minimum width=4.5pt] 
			 \draw (-3.5,0) {
			 +(-0.15,0.15) node[draw, circle,fill=gray] (w1) {}
			 +(0.15,-0.15) node[draw, circle,fill=gray] (w1p) {}
			+(45:1.2) node[draw, circle,fill=green] (w2) {}
 			+(135:1.2) node[draw, circle,fill=blue] (w3) {}
 			+(225:1.2) node[draw, circle,fill=green] (w4) {}
			+(315:1.2) node[draw, circle,fill=blue] (w5) {}
 			(w2) to (w3) to (w4) to (w5) to (w2)
			(w1) to (w2) (w1p) to (w4) (w1) to[bend right] (w3) (w1) to[bend left] (w3)  (w1p) to[bend right] (w5) (w1p) to[bend left] (w5)
			(w1)+(-0.3,-0.1) node {\small $u_1$}
			(w1p)+(0.3,0.2) node {\small $u'_1$}
 			(w2)+(0,0.3) node {\small $u_2$}
 			(w3)+(0,0.3) node {\small $u_3$}
 			(w4)+(0,-0.3) node {\small $u_4$}
			(w5)+(0,-0.3) node {\small $u_5$}
			+(-2.3,0.9) node {\small $G$}
			+(2.5,1.1) node {\small $f$}
			+(4.4,-1.3) node {\small $\check{f}$}
			};
			 \draw (-0.4,-2.6) {
			 +(0,0) node[draw, circle,fill=gray] (u1) {}
			+(45:1.2) node[draw, circle,fill=green] (u2) {}
 			+(135:1.2) node[draw, circle,fill=blue] (u3) {}
 			+(225:1.2) node[draw, circle,fill=green] (u4) {}
			+(315:1.2) node[draw, circle,fill=blue] (u5) {}
 			(u2) to (u3) to (u4) to (u5) to (u2)
			(u1) to (u2) (u1) to (u4)  (u1) to[bend right]  (u3) (u1) to[bend left]  (u3) (u1) to[bend right]  (u5) (u1) to[bend left]  (u5)
			(u1)+(0.4,0.05) node {\small $s_G$}
 			(u2)+(0,0.3) node {\small $u_2$}
 			(u3)+(0,0.3) node {\small $u_3$}
 			(u4)+(0,-0.3) node {\small $u_4$}
			(u5)+(0,-0.3) node {\small $u_5$}
			(u1)+(0,1.3) node {\small $\check{G}$}
			};
 			\draw  (4,0) {
			+(0:-2) node[draw, circle,fill=gray] (v1) {}
 			+(120:1) node[draw, circle,fill=green] (v2) {}
 			+(240:1) node[draw, circle,fill=blue] (v3) {}
 			(v1) to (v2)
			(v1) to[bend right] (v3)
			(v1) to[bend left] (v3)
			(v2) to [bend right] (v3)
			(v2) to [bend left] (v3)
 			(v1)+(0,0.3) node {\small $s_H$}
 			(v2)+(0,0.3) node {\small $v_2$}
 			(v3)+(0,-0.3) node {\small $v_3$}
			+(-0.6,1) node {\small $H$}
			};
			\draw[|->>, line width=0.3] (-2,0) to (1.5,0);
			\draw[->, gray, line width=0.4] (-3.5,-1.3) to (-1.5,-2.6);
			\draw[|->>, dashed,red, line width=0.5] (0.7,-2.6) to (2.5,-1);
\end{tikzpicture}
	\caption{\small A surjective uniform homomorphism and its induced surjective $(V(H)\setminus s_H)$-uniform  homomorphism.}
\vspace{-3mm}
\end{figure}
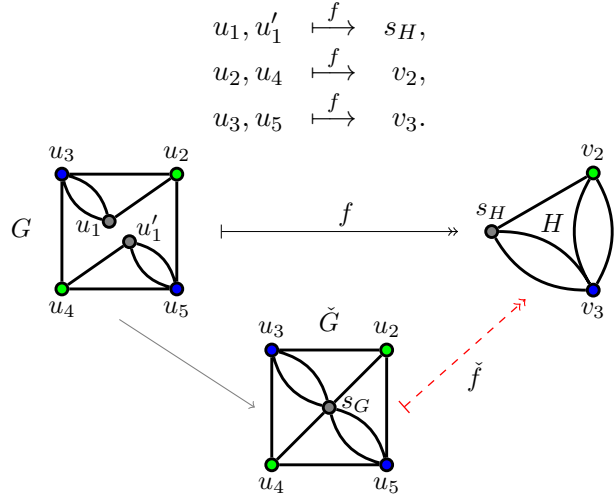

Then $\check{f}: \check{G} \rightarrow H$ is a surjective $(V(H)\setminus s_H)$-uniform  homomorphism of graphs.
\end{Remark}

\begin{Example}
In order to illustrate theorem~\ref{uniformhomeo}, consider the surjective $(V(H)\setminus s_H)$-uniform  homomorphism, 
$\check{f}: \check{G} \rightarrow H$ defined in remark~\ref{contraction}.
Using the CSandPile \footnote{CSandPile is a C++ program that computes the sandpile group of a graph. It is available by requesting to alfaromontufar@gmail.com.} 
program we can see that $SP(H,s_H)\cong \mathbb{Z}_8$ with identity ${\bf e}_{H}=(1,2)$ and generated by ${\bf c}_8=(0,3)$,
and $SP(G,s_G)\cong \mathbb{Z}_2 \oplus \mathbb{Z}_{48}$ with identity ${\bf e}_G=\widetilde{f}({\bf e}_H)=(1,2,1,2)$ and generated by 
${\bf c}_2=(2,1,2,3)$ of order two and ${\bf c}_{48}=(1,2,2,3)$ of order $48$.

\vspace{-4mm}
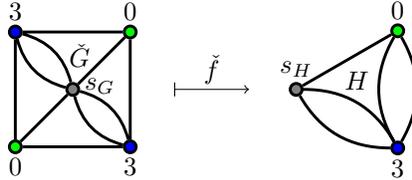
\begin{figure}[h]\centering
\begin{tikzpicture}[line width=1.1pt,scale=.9]
			\tikzstyle{every node}=[inner sep=0pt, minimum width=4.5pt] 
			 \draw (-2.1,0) {
			 +(0,0) node[draw, circle,fill=gray] (u1) {}
			+(45:1.2) node[draw, circle,fill=green] (u2) {}
 			+(135:1.2) node[draw, circle,fill=blue] (u3) {}
 			+(225:1.2) node[draw, circle,fill=green] (u4) {}
			+(315:1.2) node[draw, circle,fill=blue] (u5) {}
 			(u2) to (u3) to (u4) to (u5) to (u2)
			(u1) to (u2) (u1) to (u4)  (u1) to[bend right]  (u3) (u1) to[bend left]  (u3) (u1) to[bend right]  (u5) (u1) to[bend left]  (u5)
			(u1)+(0.4,0.05) node {\small $s_G$}
 			(u2)+(0,0.3) node {\small $0$}
 			(u3)+(0,0.3) node {\small $3$}
 			(u4)+(0,-0.3) node {\small $0$}
			(u5)+(0,-0.3) node {\small $3$}
			(u1)+(0.1,0.5) node {\small $\check{G}$}
			(u3)+(2.9,-0.55) node {\small $\check{f}$}
			};
 			\draw  (3.2,0) {
			+(0:-2) node[draw, circle,fill=gray] (v1) {}
 			+(120:1) node[draw, circle,fill=green] (v2) {}
 			+(240:1) node[draw, circle,fill=blue] (v3) {}
 			(v1) to (v2)
			(v1) to[bend right] (v3)
			(v1) to[bend left] (v3)
			(v2) to [bend right] (v3)
			(v2) to [bend left] (v3)
 			(v1)+(0,0.3) node {\small $s_H$}
 			(v2)+(0,0.3) node {\small $0$}
 			(v3)+(0,-0.3) node {\small $3$}
			+(-0.6,1) node {\small $H$}
			};
			\draw[|->, line width=0.3] (-0.6,0) to (0.5,0);
\end{tikzpicture}
\caption{\small A surjective $(V(H)\setminus s_H)$-uniform  homomorphism.}
\end{figure}
For instance, the induced mapping $\widetilde{f}: SP(H,s_H) \rightarrow SP(\check{G},s_{\check{G}})$ sends the configuration ${\bf c}_8$ 
to the configuration $\widetilde{\bf c}_8=\widetilde{f}({\bf c}_8)=(0,3,0,3)$, which generates a subgroup of order eight in $SP(\check{G},s_{\check{G}})$.
Moreover, $\widetilde{\bf c}_8={\bf c}_2\oplus 6 \cdot {\bf c}_{48}=(2,1,2,3)\oplus(2,2,2,0)$.
\end{Example}

Now, we will present a directed version of theorem~\ref{uniformhomeo}.
\begin{Theorem}\label{directeduniformhomeo}
If $G$ are a multigraph with sink $s_G$, $H$ is a multigraph with sink $s_H$, and 
$f:G\rightarrow H$ is a directed surjective $V(H)\setminus s_H$-uniform homomorphism with $f^{-1}(s_H)=\{s_G\}$.
If $\widehat{f}:\mathbb{Z}^{V(H\setminus s_H)} \rightarrow \mathbb{Z}^{V(G\setminus s_G)}$ is the mapping
given by 
\[
\widehat{f}({\bf c})_v={\bf c}_{f(v)} \text{ for all }v\in V(G)\setminus s_G,
\]
then the induced mapping $\widetilde{f}:SP(H,s_H) \rightarrow SP(G,s_G)$ 
is an injective homomorphism of groups, that is, $SP(H, s_H) \triangleleft SP(G, s_G)$.
\end{Theorem}
\begin{proof}
Clearly, $\widehat{f}$ is an an injective homomorphism of groups.
Moreover, if $L(H,s_H){\bf z}={\bf a}$, then  $L(G,s_G)\widehat{f}({\bf z})=\widehat{f}({\bf a})$.
Thus, since ${\rm det}(L(G,s_G))\neq 0$, then 
\[
\widehat{f}({\rm Im}\, L(H,s_H)) = \widehat{f}(\mathbb{Z}^{V(H\setminus s_H)})\cap {\rm Im}\, L(G,s_G).
\] 
Hence the mapping $\overline{f}:\mathbb{Z}^{V(H\setminus s_H)} \rightarrow \mathbb{Z}^{V(G\setminus s_G)}/\,{\rm Im}\, L(G,s_G)$
given by $\overline{f}({\bf a})= \widehat{f}({\bf a})\, (\text{mod }{\rm Im}\, L(G,s_G))$ has a kernel equal to ${\rm Im}\, L(H,s_H)$
and therefore the induced mapping 
\[
\widetilde{f}:SP(H,s_H)=\mathbb{Z}^{V(H\setminus s_H)}/\,{\rm Im}\, L(H,s_H) \rightarrow \mathbb{Z}^{V(G\setminus s_G)}/ \,{\rm Im}\, L(G,s_G)=SP(G,s_G)
\]
is an injective homomorphism of groups.
\end{proof}

\begin{Remark}
Note that, if $f:G\rightarrow H$ is a directed surjective $V(H) \setminus s_H$-uniform homomorphism with $f^{-1}(s_H)=\{s_G\}$ 
and ${\bf a}$ is an eigenvector of $L(H, s_H)$ for the eigenvalue $\lambda$,
then $\widehat{f}({\bf a})$ is an eigenvector of $L(G, s_G)$ for $\lambda$.

Also, note that in the directed case the mapping $\widetilde{f}$ defined in theorem~\ref{directeduniformhomeo} is not a natural homomorphism of sandpile groups
in the sense that it does not necessarily send recurrent configurations to recurrent configurations.
\end{Remark}

The next corollary is an application of theorem~\ref{directeduniformhomeo}.

\begin{Corollary}\label{regular}
If $B_{r,t}$ is a bipartite graph with bipartition $V=V_1\cup V_2$ and
\[
d(u)=
\begin{cases}r & \text{ if } u\in V_1,\\
t & \text{ if } u\in V_2,
\end{cases}
\]
then ${\mathbb Z}_{r+t+1}  \triangleleft SP(c(B_{r,t}),s_{B_{r,t}})$.
\end{Corollary}
\begin{proof}
Let ${\mathcal K}_2(r,t)$ be the multigraph with $V=\{v_1,v_2\}$ as set of vertices and $m_{v_1,v_2}=r$, $m_{v_2,v_1}=t$.
Let $f:c(B_{r,t})\rightarrow c({\mathcal K}_2(r,t))$ be the mapping given by

\[
f(v)=
\begin{cases}
v_1 & \text { if } v\in V_1,\\
v_2 & \text { if } v\in V_2,\\
s_{{\mathcal K}_2(r,t)} & \text{ if } v=s_{B_{r,t}}.
\end{cases}
\]
Since $f$ is a surjective $\{v_1,v_2\}$-uniform homomorphism, then by theorems \ref{directeduniformhomeo} and \ref{generators} we have that
\[
{\mathbb Z}_{r+t+1}\cong SP(c({\mathcal K}_2(r,t)))  \triangleleft SP(c(B_{r,t}),s_{B_{r,t}}).
\]
\end{proof}

\begin{Remark}\label{generalized}
A \textit{weak $V$-uniform homomorphism} is a mapping $f:V(G)\rightarrow V(H)$ such that for all $x\in V$ and $y\in V(H)$ with $x\neq y$
\[
d_{G[\{u\}\cup S_y]}(u)=m_{x,y} \text{ for all } u \in S_x
\]
(that is, the sets $S_x$ are not necessarily stable) and $f: V(G)\setminus f^{-1}(V)\rightarrow V(H)\setminus V$ is the identity isomorphism.
In this case the induced mapping $\widetilde{f}$ does not send recurrent configurations to recurrent configurations, 
but the mapping $\hat{f}({\bf c})=[\widetilde{f}({\bf c})]$
(where $[\widetilde{f}({\bf c})]$ is the unique recurrent configuration of $G$ such that 
$s(\widetilde{f}({\bf c})+r)=[\widetilde{f}({\bf c})]$ for some non-zero configuration $r$)
is an injective homomorphism of groups.
\end{Remark}

\begin{Remark}\label{bicycles}
The group of bicycles of a graph $G$ over an abelian group $A$, denoted by $B(G,A)$, 
consists of the edge weightings of $G$ over $A$ that are both cycles and cocycles of $G$ and the entry by entry sum.
The group of bicycles and the sandpile group of a graph are closely related.
For instance, $B(G, A)={\rm Hom}_{\mathbb Z}(SP(G),A)$.
Moreover, if either $A=\mathbb{Q}/ \mathbb{Z}$ or $A=\mathbb{Z}_{|SP(G)|}$, 
then the group $B(G, A)$ of bicycles of $G$ is isomorphic to $SP(G)$.  

Let $G$, $H$ be connected multigraphs and $V(H)=\{u_1,\cdots, u_{|V(H)|}\}$ be the vertex set of $H$. 
We say that $G$ is \textit{divisible} by $H$ (see~\cite[page 9]{bicycles}) if the vertices of $G$ 
can be partitioned into $|V(H)|$ classes $U_1,\cdots, U_{|V(H)|}$,
such that for $1\leq i, j\leq |V(H)|$ a vertex $v$ in $U_i$ is either joined only to vertices of 
$U_i$ or for every $i\neq j$ is joined to exactly $m_{u_i,u_j}$ vertices of $U_{j}$ (and any number of vertices in $U_i$).

Note that the concepts of divisibility and weak $V(G)$-uniform homomorphism are closed related.
Clearly, if $G$ is divisible by $H$, then there exists a weak $V(H)$-uniform homomorphism $f$
between $G$ and $H$.
However, if $\check{G}$ and $H$ are the graphs defined in remark~\ref{contraction}, then $\check{G}$ is not divisible by $H$
but there exists a surjective $(V(H)\setminus s_H)$-uniform homomorphism $\check{f}$ between $\check{G}$ and $H$.
Also, it is not difficult to see that the cycle with four vertices with an added pendant edge 
(that is, $E(G)=\{x_1x_2,x_2x_3,x_3x_4,x_4x_1,x_1x_5\}$) is divisible by $\mathcal{K}_2(2)$ but
there not exists a uniform homomorphism between them.

Theorem 5.7 in~\cite{bicycles} says that if $G$ is divisible by $H$, then $B(H,\mathbb{Z}_k)$ is a subgroup of $B(G,\mathbb{Z}_k)$ for all $k\in \mathbb{Z}$.
That is, Theorems \ref{uniformhomeo}, \ref{directeduniformhomeo}, and ~\cite[Theorem 5.7]{bicycles} shows injections between groups
induced by some class of morphism between graphs.
\end{Remark} 

%==================================================================% 
%=======================The Cartesian Product of graphs=================% 
%==================================================================% 

\section{The sandpile group of the cartesian product of graphs}

The sandpile group of the cartesian product of graphs has been studied by several authors,
see for instance~\cite{Bai,p4cn,k3cn,cartesian,kmpn}. 
In this section we define the cartesian product of configurations and we prove that the cartesian product
of recurrent configurations is a recurrent configuration.
After that, we prove that: if ${\bf e}_H\in SP(c(H), s_{c(H)})$ is the identity of the sandpile group of the cone of $H$,
then the mapping $\widetilde{\pi}_G: SP(c(G), s_{c(G)})\rightarrow SP(c(G\Box H), s_{c(G\Box H)})$
given by
\[
\widetilde{\pi}_G ({\bf a})={\bf a}\Box {\bf e}_H,
\]
is an injective homomorphism of groups.

The \textit{cartesian product} of $G$ and $H$, denoted by $G\Box H$ is the graph 
with $V(G) \times V(H)$ as its vertex set and two vertices $u_1 v_1$ and $u_2 v_2$ are adjacent in $G \Box H$ 
if and only if either $u_1 = u_2$ and $v_1 v_2\in E(H)$, or $v_1 = v_2$ and $u_1 u_2\in E(G)$.
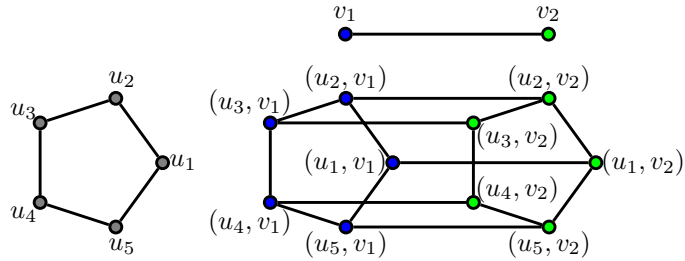
\begin{figure}[h]\centering
		\begin{tikzpicture}[line width=1.1,scale=0.9]
		\tikzstyle{every node}=[inner sep=0pt, minimum width=4.5pt] 
 		\draw (-3.4,0) {
		+(0:1) node[draw, circle, fill=gray] (u1) {}
 		+(72:1) node[draw, circle, fill=gray] (u2) {}
 		+(144:1) node[draw, circle, fill=gray] (u3) {}
 		+(216:1) node[draw, circle, fill=gray] (u4) {}
 		+(288:1) node[draw, circle, fill=gray] (u5) {}
		+(0:1.3) node {\small $u_1$}
 		+(72:1.3) node {\small $u_2$}
 		+(144:1.3) node {\small $u_3$}
 		+(216:1.3) node {\small $u_4$}
 		+(288:1.3) node {\small $u_5$}
 		(u1) to (u2) to (u3) to (u4) to (u5) to (u1)
		};
		\draw (0,1.9) {
		+(0.3,0) node[draw, circle, fill=blue] (v1) {}
 		+(3.3,0) node[draw, circle, fill=green] (v2) {}
		(v1)+(0,0.3) node {\small $v_1$}
 		(v2)+(0,0.3) node {\small $v_2$}
 		(v1) to (v2)
		};
		\draw (0,0) {
		+(0:1) node[draw, circle, fill=blue] (v1) {}
 		+(72:1) node[draw, circle, fill=blue] (v2) {}
 		+(144:1) node[draw, circle, fill=blue] (v3) {}
 		+(216:1) node[draw, circle, fill=blue] (v4) {}
 		+(288:1) node[draw, circle, fill=blue] (v5) {}
		(v1)+(3,0) node[draw, circle, fill=green] (v11) {}
 		(v2)+(3,0) node[draw, circle, fill=green] (v12) {}
 		(v3)+(3,0) node[draw, circle, fill=green] (v13) {}
 		(v4)+(3,0) node[draw, circle, fill=green] (v14) {}
 		(v5)+(3,0) node[draw, circle, fill=green] (v15) {}
		(v1)+(-0.7,0) node {\small $(u_1,v_1)$}
 		(v2)+(0,0.3) node {\small $(u_2,v_1)$}
 		(v3)+(-0.3,0.3) node {\small $(u_3,v_1)$}
 		(v4)+(-0.3,-0.3) node {\small $(u_4,v_1)$}
 		(v5)+(0,-0.25) node {\small $(u_5,v_1)$}
		(v11)+(0.7,0) node {\small $(u_1,v_2)$}
 		(v12)+(0,0.3) node {\small $(u_2,v_2)$}
 		(v13)+(0.65,-0.2) node {\small $(u_3,v_2)$}
 		(v14)+(0.65,0.2) node {\small $(u_4,v_2)$}
 		(v15)+(0,-0.25) node {\small $(u_5,v_2)$}
 		(v1) to (v2) to (v3) to (v4) to (v5) to (v1)
		(v11) to (v12) to (v13) to (v14) to (v15) to (v11)
		(v1) to (v11) (v2) to (v12) (v3) to (v13) (v4) to (v14) (v5) to (v15)
		};
		\end{tikzpicture}
\caption{\small Cartesian product of $C_5$ and $\mathcal{K}_2$.} 
\end{figure}

Let $\pi_{G}: G\Box H \rightarrow G$ and $\pi_{H}:G\Box H \rightarrow H$ be the projection mappings, given by
\[
\pi_{G}(u,v)=u \text{ for all } (u,v)\in V(G\Box H) \text{ and } \pi_{H}(u,v)=v \text{ for all } (u,v)\in V(G\Box H).
\]

Thus, it is not difficult to see that the mappings $\pi_{G}$ and $\pi_{H}$ are weak surjective uniform homomorphisms of graphs. 
For the rest of this section, let $s_{c(G)}\in V(c(G))\setminus V(G)$,  $s_{c(H)}\in V(c(H))\setminus V(H)$ and $s_{c(G\Box H)}\in V(c(G \Box H))\setminus V(G \Box H)$.

\medskip

Now, let ${\bf a}\in \mathbb{N}^{V(G)}$ be a configuration of $c(G)$, ${\bf b}\in \mathbb{N}^{V(H)}$ be a configuration of $c(H)$, 
\vspace{-5mm}
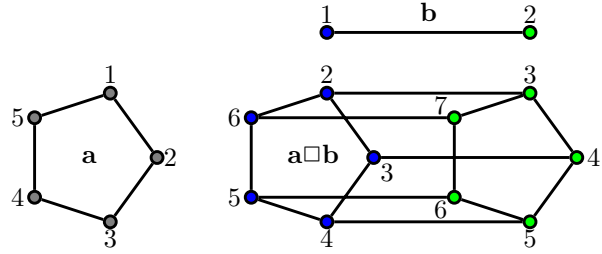
\begin{figure}[h] \centering
	\begin{tikzpicture}[line width=1.1,scale=0.9]
		\tikzstyle{every node}=[inner sep=0pt, minimum width=4.5pt] 
		\draw (0,0) {
		+(0:1) node[draw, circle, fill=blue] (v1) {}
		+(72:1) node[draw, circle, fill=blue] (v2) {}
		+(144:1) node[draw, circle, fill=blue] (v3) {}
		+(216:1) node[draw, circle, fill=blue] (v4) {}
		+(288:1) node[draw, circle, fill=blue] (v5) {}
		(v1)+(3,0) node[draw, circle, fill=green] (v11) {}
		(v2)+(3,0) node[draw, circle, fill=green] (v12) {}
		(v3)+(3,0) node[draw, circle, fill=green] (v13) {}
		(v4)+(3,0) node[draw, circle, fill=green] (v14) {}
		(v5)+(3,0) node[draw, circle, fill=green] (v15) {}
 		(v1) to (v2) to (v3) to (v4) to (v5) to (v1)
		(v11) to (v12) to (v13) to (v14) to (v15) to (v11) (v1) to (v11) (v2) to (v12) (v3) to (v13) (v4) to (v14) (v5) to (v15)
		(v1)+(0.20,-0.22) node {\small $3$}
		(v2)+(0.0,0.27) node {\small $2$}
		(v3)+(-0.25,0.0) node {\small $6$}
		(v4)+(-0.25,0.0) node {\small $5$}
		(v5)+(0.0,-0.25) node {\small $4$}
		(v11)+(0.25,0) node {\small $4$}
		(v12)+(0.0,0.27) node {\small $3$}
		(v13)+(-0.20,0.18) node {\small $7$}
		(v14)+(-0.2,-0.18) node {\small $6$}
		(v15)+(0.0,-0.25) node {\small $5$}
		+(-3.2,1) node {\small ${\bf a}\Box{\bf b}$}
		};
		\draw (-3.2,0) {
		+(v1) node[draw, circle, fill=gray] (v16) {}
		+(v2)node[draw, circle, fill=gray] (v17) {}
		+(v3) node[draw, circle, fill=gray] (v18) {}
		+(v4) node[draw, circle, fill=gray] (v19) {}
		+(v5) node[draw, circle, fill=gray] (v20) {}
		(v16) -- (v17) -- (v18) -- (v19) -- (v20) -- (v16)
		+(-1,0) node {\small ${\bf a}$}
		(v16)+(0.2,0) node {\small $2$}
		(v17)+(0.0,0.27) node {\small $1$}
		(v18)+(-0.25,0.0) node {\small $5$}
		(v19)+(-0.25,0.0) node {\small $4$}
		(v20)+(0.0,-0.25) node {\small $3$}
		};		
		\draw (0,2.8) {
		+(v5)node[draw, circle, fill=blue] (v21) {}
		+(v15) node[draw, circle, fill=green] (v22) {}
		(v21) to (v22)
		(v21)+(0,0.27) node {\small $1$}
		(v22)+(0,0.27) node {\small $2$}
		+(-1.5,0.3) node {\small ${\bf b}$}
		};
	\end{tikzpicture}
	\caption{\small Cartesian product of configurations.}
	\label{cartproconfi}
\end{figure}
and let ${\bf a}\Box {\bf b}\in \mathbb{N}^{V(G\Box H)}$ be the configuration of $c(G\Box H)$ given by 
\[
({\bf a}\Box {\bf b})_{(u,v)}={\bf a}_u+{\bf b}_v  \text{ for all } u\in V(G) \text{ and }  v \in V(H).
\]

The following lemma shows that the cartesian product of configurations of $c(G)$ and $c(H)$ 
is compatible with the toppling operators of $c(G)$, $c(H)$ and $c(G\Box H)$.

\begin{Lemma}\label{recurrent}
Let $G$ and $H$ be multigraphs, ${\bf a}\in \mathbb{N}^{V(G)}$ be a configuration of $c(G)$, and ${\bf b}\in \mathbb{N}^{V(H)}$ be a configuration of $c(H)$. 
Then
\begin{description}
\item[$(i)$] If ${\bf a}$ and ${\bf b}$ are stable configurations, then ${\bf a}\Box {\bf b}$ is a stable configuration of $c(G\Box H)$,
\item[$(ii)$] If ${\bf a}$ and ${\bf b}$ are recurrent configurations, then ${\bf a}\Box {\bf b}$ is a recurrent configuration of $c(G\Box H)$.
\end{description} 
\end{Lemma}
\begin{proof}
$(i)$
If ${\bf a}$ and ${\bf b}$ are stable configurations of $c(G)$ and $c(H)$ respectively, 
then 
\[
{\bf a}_u \leq deg_{c(G)}(u)-1 \text{ for all } u\in V(G) \text{ and }{\bf b}_v \leq deg_{c(H)}(v)-1 \text{ for all } v\in V(H).
\]
Hence ${\bf a}\Box {\bf b}_{(u,v)}={\bf a}_u+{\bf b}_v \leq deg_{c(G)}(u)+deg_{c(H)}(v)-2=deg_{c(G\Box H)}((u,v))-1$, 
that is, ${\bf a}\Box {\bf b}$ is a stable configuration of $c(G\Box H)$.

$(ii)$
We will use the burning algorithm~\ref{burning} to prove the second part of this lemma.
Since, the sink $s_G$ of $c(G)$ is adjacent to all the vertices of $G$, then $\sum_{i=1}^n \Delta_{u_i}={\bf 1}$.

\begin{Claim}
Let ${\bf a}$ be a recurrent configuration of $c(G)$ and ${\bf b}$ be a recurrent configuration $c(H)$.
Also, let
\[
{\bf a}_{i}=
\begin{cases}
{\bf a}+{\bf 1} &\text{ if } i=1\\ 
{\bf a}_{i-1}-\Delta_{u_{i-1}} &\text{ if } i=2,\ldots,n, 
\end{cases}
\text{ and } 
{\bf b}_{i}=
\begin{cases}
{\bf b}+{\bf 1} &\text{ if } j=1\\ 
{\bf b}_{i-1}-\Delta_{v_{i-1}} & \text{ if } j=2,\ldots,m,
\end{cases}
\]   
such that the vertex $u_{i}$ is an unstable vertex in ${\bf a}_{i}$ for all $i=1,\ldots,n$  and the vertex $v_{j}$ is an unstable vertex in ${\bf b}_{j}$ for all $j=1,\ldots,m$.
If ${\bf c}={\bf a}\Box {\bf b}$, ${\bf c}_{(1,1)}={\bf a}\Box {\bf b}+{\bf 1}={\bf a}_1\Box {\bf b}={\bf a}\Box {\bf b}_1$, and
\[
{\bf c}_{(i,j)}=
\begin{cases}
{\bf c}_{(i-1,m)}-\Delta_{(u_{i-1},v_m)} & \text{ if } i=2,\ldots,n \text{ and } j=1,\\
{\bf c}_{(i,j-1)}-\Delta_{(u_i,v_{j-1})} & \text{ otherwise, }
\end{cases}
\] 
then the vertex $(u_i,v_j)$ is an unstable vertex in ${\bf c}_{(i,j)}$ for all $i=1,\ldots,n$ and $j=1,\ldots,m$.
\end{Claim}
\begin{proof}
Since the vertex $u_{i}$ is an unstable vertex in ${\bf a}_{i}$ for all $i=1,\ldots,n$ and the vertex $v_{j}$ is an unstable vertex in ${\bf b}_{j}$ for all $j=1,\ldots,m$,
then $({\bf a}_{i})_{u_{i}} \geq  {\rm deg}_{c(G)} (u_i)$ for all $i=1,\ldots,n$ and $({\bf b}_{j})_{v_{j}} \geq  {\rm deg}_{c(H)} (v_j)$ for all $j=1,\ldots,m$.

\medskip

Now, ${\bf c}_{(i,1)}={\bf c}_{(1,1)}-\sum_{1\le k\leq i-1} \sum_{1\le l\leq m}  \Delta_{(u_{k},v_l)}
=({\bf a}_1- \sum_{1\le k\leq i-1}  \Delta_{u_{k}})\Box {\bf b}={\bf a}_{i}\Box {\bf b}$ for all $i=1,\ldots,n$.
Thus, 
\begin{eqnarray*}
({\bf c}_{(i,1)})_{(u_i,v_1)}&=&({\bf a}_{i}\Box {\bf b})_{(u_i,v_1)}=({\bf a}_{i})_{u_i}+{\bf b}_{v_1}=({\bf a}_{i})_{u_i}+({\bf b}_1)_{v_1}-1
\geq {\rm deg}_{c(G)} (u_i)+{\rm deg}_{c(H)} (v_1)-1\\ &=&{\rm deg}_{c(G\Box H)} ((u_i,v_1)) \text{ for all }i=1,\ldots,n.
\end{eqnarray*}
Moreover, since  ${\bf c}_{(i,j)}={\bf a}_{i}\Box {\bf b}-\sum_{1\le l\leq j}  \Delta_{(u_{i},v_l)}=
({\bf a}_{i}-{\bf 1})\Box {\bf b}_1-\sum_{1\le l\leq j}  \Delta_{(u_{i},v_l)}$, then   
\[
({\bf c}_{(i,j)})_{(u_i,v_j)}=({\bf a}_{i})_{u_i}+ ({\bf b}_j)_{v_j}-1
\geq {\rm deg}_{c(G)} (u_i)+{\rm deg}_{c(H)} (v_j)-1={\rm deg}_{c(G\Box H)} ((u_i,v_j))
\]
for all $i=1,\ldots,n$ and $j=1,\ldots,m$.

Therefore $(u_i,v_j)$ is an unstable vertex of ${\bf c}_{(i,j)}$ for all $i=1,\ldots,n$ and $j=1,\ldots,m$.
\end{proof}
Finally, by using part $(i)$ of this lemma and the previous claim we obtain that ${\bf a}\Box {\bf b}$ is recurrent.
\end{proof}

The next example is useful to illustrate the previous theorem:

\begin{Example}
Let $G\cong H\cong \mathcal{K}_2$ with $V(G)=\{u_1,u_2\}$ and $V(H)=\{v_1,v_2\}$ as vertex sets,
${\bf a}=(1,1)$ be a recurrent configuration of $c(G)$ and ${\bf b}=(1,0)$ be a recurrent configuration of $c(H)$.
 
Hence ${\bf c}=(2,1,2,1)=(1,1)\Box (1,0)$ is a recurrent configuration of $c(G\Box H)$, as is shown in figure~$7$.
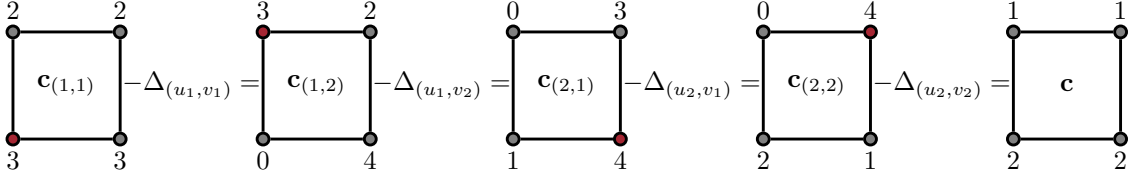
\begin{figure}[h] \centering
\label{toppling}
	\begin{tikzpicture}[line width=1.1pt,scale=0.95]
		\tikzstyle{every node}=[inner sep=0pt, minimum width=4.5pt] 
		\draw (-7.0,0) {
		+(0,0) node[draw, circle, fill=Magenta] (v1) {}
		+(1.5,0) node[draw, circle, fill=gray] (v2) {}
		+(1.5,1.5) node[draw, circle, fill=gray] (v3) {}
		+(0,1.5) node[draw, circle, fill=gray] (v4) {}
		(v1) -- (v2) -- (v3)--(v4)--(v1)
		(v1)+(0,-0.3) node {\small $3$}
		(v2)+(0,-0.3) node {\small $3$}
		(v3)+(0,0.3) node {\small $2$}
		(v4)+(0,0.3) node {\small $2$}
		+(0.75,-0.75) node {\small ${\bf c}_{(1,1)}$}
		+(2.5,-0.75) node {\small $-\Delta_{(u_1,v_1)}=$}
		}; 

		\draw (-3.5,0) {
		+(0,0) node[draw, circle, fill=gray] (v5) {}
		+(1.5,0) node[draw, circle, fill=gray] (v6) {}
		+(1.5,1.5) node[draw, circle, fill=gray] (v7) {}
		+(0,1.5) node[draw, circle, fill=Magenta] (v8) {}
		(v5) -- (v6) -- (v7)--(v8)--(v5)
		(v5)+(0,-0.3) node {\small $0$}
		(v6)+(0,-0.3) node {\small $4$}
		(v7)+(0,0.3) node {\small $2$}
		(v8)+(0,0.3) node {\small $3$}
		+(0.75,-0.75) node {\small ${\bf c}_{(1,2)}$}
		+(2.5,-0.75) node {\small $-\Delta_{(u_1,v_2)}=$}
		}; 	
			
		\draw (0,0) {
		+(0,0) node[draw, circle, fill=gray] (v9) {}
		+(1.5,0) node[draw, circle, fill=Magenta] (v10) {}
		+(1.5,1.5) node[draw, circle, fill=gray] (v11) {}
		+(0,1.5) node[draw, circle, fill=gray] (v12) {}
		(v9) -- (v10) -- (v11)--(v12)--(v9)
		(v9)+(0,-0.3) node {\small $1$}
		(v10)+(0,-0.3) node {\small $4$}
		(v11)+(0,0.3) node {\small $3$}
		(v12)+(0,0.3) node {\small $0$}
		+(0.75,-0.75) node {\small ${\bf c}_{(2,1)}$}
		+(2.5,-0.75) node {\small $-\Delta_{(u_2,v_1)}=$}
		};
			
 		\draw (3.5,0) {
		+(0,0) node[draw, circle, fill=gray] (v13) {}
		+(1.5,0) node[draw, circle, fill=gray] (v14) {}
		+(1.5,1.5) node[draw, circle, fill=Magenta] (v15) {}
		+(0,1.5) node[draw, circle, fill=gray] (v16) {}
		(v13) -- (v14) -- (v15)--(v16)--(v13)
		(v13)+(0,-0.3) node {\small $2$}
		(v14)+(0,-0.3) node {\small $1$}
		(v15)+(0,0.3) node {\small $4$}
		(v16)+(0,0.3) node {\small $0$}
		+(0.75,-0.75) node {\small ${\bf c}_{(2,2)}$}
		+(2.5,-0.75) node {\small $-\Delta_{(u_2,v_2)}=$}
		};
		
		\draw (7.0,0) {
		+(0,0) node[draw, circle, fill=gray] (v17) {}
		+(1.5,0) node[draw, circle, fill=gray] (v18) {}
		+(1.5,1.5) node[draw, circle, fill=gray] (v19) {}
		+(0,1.5) node[draw, circle, fill=gray] (v20) {}
		(v17) -- (v18) -- (v19)--(v20)--(v17)
		(v17)+(0,-0.3) node {\small $2$}
		(v18)+(0,-0.3) node {\small $2$}
		(v19)+(0,0.3) node {\small $1$}
		(v20)+(0,0.3) node {\small $1$}
		+(0.75,-0.75) node {\small ${\bf c}$}
		};					
\end{tikzpicture}
\caption{\small The topplings of the configuration ${\bf c}=(3,2,3,2)$ of $c(C_4)$.}
\end{figure}
\vspace{-3mm}
\end{Example}

The next theorem shows that the mappings $\pi_{G}$ and $\pi_{H}$ induce homomorphisms of groups 
between the sandpile groups of the cones of $G$ and $G\Box H$, and $H$ and $G\Box H$; respectively.

\begin{Theorem}\label{cartesian}
Let $G$ and $H$ be two multigraphs, and ${\bf e}_H$ be the identity of the sandpile group of the cone of $H$.
Then the mapping $\widetilde{\pi}_G: SP(c(G), s_{c(G)})\rightarrow SP(c(G\Box H), s_{c(G\Box H)})$
given by
\[
\widetilde{\pi}_G ({\bf a})={\bf a}\Box {\bf e}_H,
\]
is an injective homomorphism of groups.
\end{Theorem}
\begin{proof}
Since ${\bf e}_H$ is recurrent, then using lemma~\ref{recurrent} $(ii)$,
$\widetilde{\pi}_G ({\bf a})={\bf a}\Box {\bf e}_H$ is a recurrent configuration of $c(G\Box H)$ for all ${\bf a}\in SP(c(G), s_{c(G)})$;
that is, the mapping $\widetilde{\pi}_G$ is well defined. 

Now, we will prove that $\widetilde{\pi}_G$ is a homomorphism of groups. 
Let ${\bf a},{\bf b}\in SP(c(G), s_{c(G)})$, then
\begin{eqnarray*}
\widetilde{\pi}_G ({\bf a} \oplus {\bf b})&=&({\bf a}\oplus {\bf b})\Box {\bf e}_H=
s({\bf a}+ {\bf b})\Box {\bf e}_H=({\bf a}+ {\bf b})\Box {\bf e}_H \,(\text{mod }L(c(G\Box H),s_{c(G\Box H)}))\\
&=&{\bf a}\Box {\bf e}_H+{\bf b}\Box {\bf e}_H = s({\bf a}\Box {\bf e}_H + {\bf b}\Box{\bf e}_H) \,(\text{mod }L(c(G\Box H),s_{c(G\Box H)}))=
{\bf a}\Box {\bf e}_H \oplus {\bf b}\Box{\bf e}_H\\
&=&\widetilde{\pi}_G ({\bf a})\oplus \widetilde{\pi}_G ({\bf b}),
\end{eqnarray*}
and therefore $\widetilde{\pi}_G$ is a homomorphism of groups.
	
Finally, $\widetilde{\pi}_G ({\bf a})=\widetilde{\pi}_G ({\bf b})$ 
if and only if ${\bf a}\Box {\bf e}_H={\bf b}\Box {\bf e}_H$ if and only if ${\bf a}={\bf b}$,
and therefore $\widetilde{\pi}_G$ is an injective homomorphism of groups.
\end{proof}

\begin{Example}
Using the CSandPile program we get that, $SP(c(\mathcal{K}_2), s_{c(\mathcal{K}_2)})=\mathbb{Z}_3$ is generated by $(1,0)$ with identity $(1,1)$,
$SP(c(C_5),s_{c(C_5)})=\mathbb{Z}_{11}^{2}$ is generated by $(2,1,1,1,1)$ and $(1,2,1,1,1)$ with identity $e=(2,2,2,2,2)$ (Also see~\cite[page 5]{corirossin}), 
and $SP(c(C_5\Box \mathcal{K}_2))=\mathbb{Z}_{11\cdot 29} \oplus \mathbb{Z}_{3\cdot 11\cdot 29}$. 

Moreover, using the mapping $\pi_{\mathcal{K}_2}$ we have that
\[
\pi_{\mathcal{K}_2}(1,0)=(3,3,3,3,3,2,2,2,2,2)
\]
is a generator of a subgroup of $SP(c(C_5\Box \mathcal{K}_2))$ isomorphic to $\mathbb{Z}_3$,
and using the mapping $i_{C_5}$ we have that
\[
\pi_{C_5}(2,1,1,1,1)=(3,2,2,2,2,3,2,2,2,2) \text{ and }  \pi_{C_5}(1,2,1,1,1)=(2,3,2,2,2,2,3,2,2,2)
\]
are generators of subgroups of $SP(c(C_5\Box \mathcal{K}_2))$ isomorphic to $\mathbb{Z}_{11}$.
\end{Example}

\begin{Remark}
If $n> 1$ and ${\bf e}_H$ is the identity of $SP(c_n(G), s_{c_n(G)})$, then the mapping given by
\[
\widetilde{\pi}_G ({\bf a})={\bf a}\Box {\bf e}_H
\]
does not necessarily send stable configurations to stable configurations.
For instance, the vector $(3,3)$ is the identity of $c_3(Q_1)$ and $(6,6,6,6)=\widetilde{\pi}_G ((3,3))=(3,3)\Box (3,3)$ 
is a non stable configuration of $c_3(Q_2)$.
However, the non canonical mapping $$\hat{\pi}_G: SP(c_n(G), s_{c_n(G)})\rightarrow SP(c_n(G\Box H), s_{c_n(G\Box H)})$$
given by $\hat{\pi}_G ({\bf a})=[\widetilde{\pi}_G ({\bf a})]$ is an injective homomorphism of groups.
\end{Remark}

%==================================================================% 
%=============================Hipercubo=========================% 
%==================================================================% 

\section{The sandpile group of $c(Q_d)$}

The hypercube of dimension $d$ is the cartesian product of $d$ copies of the complete graph with two vertices ${\mathcal K}_2$.
The structure of the sandpile group of the hypercube is complex, see for intance~\cite{Bai} for a description of the Sylow $p$-group of $SP(Q_d)$ 
when $p$ is odd and ~\cite{cartesian} for a description of the cartesian product of complete graphs in general.

In this section we give an explicit combinatorial and algebraic description of a set of generators of the sandpile 
group of the cone of the hypercube of dimension $d$, see theorem~\ref{cQd}.
We will use mainly theorems~\ref{uniformhomeo} and ~\ref{cartesian}, developed in previous sections, to get 
a description of the sandpile group of the cone of the hypercube of dimension $d$.

\medskip

First of all, we will fix some notation that will be needed in order to establish the main theorem, let
\[
Q_d=\Box_{i=1}^d {\mathcal K}_2=\underbrace{{\mathcal K}_2 \Box \cdots \Box {\mathcal K}_2}_{d- \text{copies of } {\mathcal K}_2},
\]
with vertex set $V(Q_d)=\{v_{\bf a} \, | \, {\bf a}\in \{0,1\}^d\}$ and edge set
\[
E(Q_d)=\{v_{\bf a}v_{{\bf a}'} \, | \, {\bf a},{\bf a}'\in \{0,1\}^d \text{ and }  {\bf a}-{\bf a}'=\pm e_i \text{ for some } 1\leq i \leq d \}.
\]
Moreover, for all $\beta \in \{0,1\}^d$, let $Q_{\beta}=Q_d[\{ v_{\bf a} \, | \, {\rm supp}({\bf a}) \subseteq {\rm supp}(\beta) \}]$,
be an induced subgraph of $Q_d$, where ${\rm supp}(c)=\{i \mid c_i\neq 0 \}$. 
It is not difficult to note that: 
\begin{itemize}
\item $Q_{\beta} = \Box_{i\in {\rm supp}(\beta)} Q_{e_i}\cong \Box_{i=1}^{|\beta|} {\mathcal K}_2=Q_{|\beta|} $, where $|\beta|=\sum_{i=1}^d \beta_i$,
\item $Q_{\beta'}=Q_{\beta} \Box Q_{\beta'-\beta}$ for all ${\rm supp}(\beta) \subset {\rm supp}(\beta')$, in particular
$Q_d=Q_{(1,\ldots,1)}=Q_\beta \Box Q_{{\bf 1}-\beta}$ for all $\beta \in \{0,1\}^d$.
\end{itemize}

\vspace{-6mm}
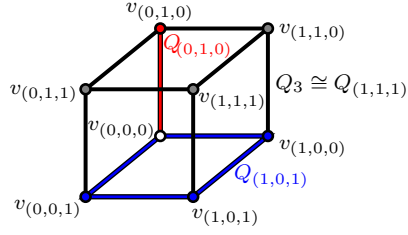
\begin{figure}[h] \centering
\begin{tikzpicture}[line width=1pt, scale=1.1]
\path (0,0) coordinate (v1);
\path (1.3,0) coordinate (v2);
\path (226:1.3)+(1.3,0.2) coordinate (v3);
\path (226:1.3)+(0,0.2) coordinate (v4);

\path (0,0)+(0,1.3) coordinate (v5);
\path (1.3,0)+(0,1.3) coordinate (v6);
\path (226:1.3)+(1.3,1.5) coordinate (v7);
\path (226:1.3)+(0,1.5) coordinate (v8);

\draw[line width=2pt] (v1)--(v2)--(v3)--(v4)--(v1);
\draw[blue,line width=1pt] (v1)--(v2)--(v3)--(v4)--(v1);

\draw[line width=2pt] (v1)--(v5);
\draw[red,line width=1pt] (v1)--(v5);

\draw[line width=1.5pt] (v5)--(v6)--(v7)--(v8)--(v5);

\draw[line width=1.5pt] (v2)--(v6) (v3)--(v7) (v4)--(v8);

\draw[fill=white] (v1) circle(1.7pt);
\draw[fill=blue] (v2) circle(1.7pt) (v3) circle(1.7pt) (v4) circle(1.7pt);
\draw[fill=red] (v5) circle(1.7pt);
\draw[fill=gray] (v6) circle(1.7pt) (v7) circle(1.7pt) (v8) circle(1.7pt);

\draw (v1)+(-0.46,0.07) node {\scriptsize $v_{(0,0,0)}$};
\draw (v2)+(0.55,-0.15) node {\scriptsize $v_{(1,0,0)}$};
\draw (v3)+(0.40,-0.25) node {\scriptsize $v_{(1,0,1)}$};
\draw (v4)+(-0.45,-0.15) node {\scriptsize $v_{(0,0,1)}$};
\draw (v5)+(-0.0,0.20) node {\scriptsize $v_{(0,1,0)}$};
\draw (v6)+(0.55,-0.05) node {\scriptsize $v_{(1,1,0)}$};
\draw (v7)+(0.48,-0.15) node {\scriptsize $v_{(1,1,1)}$};
\draw (v8)+(-0.50,-0.05) node {\scriptsize $v_{(0,1,1)}$};
\draw[text=red]  (v3)+(0.06,1.82) node {\scriptsize $Q_{\mbox{} \hspace{-0.7mm} (0,1,0)}$};
\draw[text=blue]  (v3)+(0.95,0.25) node {\scriptsize $Q_{(1,0,1)}$};
\draw (v2)+(0.9,0.6) node {\scriptsize $Q_3\cong Q_{(1,1,1)}$};
\end{tikzpicture}
\caption{\small The hypercube $Q_3\cong Q_{(1,1,1)}$, where the hypercube \textcolor{blue}{$Q_{(1,0,1)}\cong Q_2$} is colored in blue, 
and the hypercube \textcolor{red}{$Q_{(0,1,0)}\cong Q_1$} is colored in red.}
\end{figure}

Now, for all ${\bf 0}\neq \beta\in \{0,1\}^d$, let $f_{\beta}:c(Q_{\beta}) \rightarrow c(\mathcal{K}_2(|\beta|))$ 
be the surjective $V(\mathcal{K}_2(|\beta|))$-uniform homomorphism of graphs given by 
\[
f_{\beta}(v)=
\begin{cases}
v_1 & \text{ if } v=v_{\bf a} \text{ and } \beta\cdot {\bf a} \text{ is even},\\
v_2 & \text{ if } v=v_{\bf a} \text{ and } \beta\cdot {\bf a} \text{ is odd},\\
s_{c(\mathcal{K}_2(|\beta|))}& \text{ if } v=s_{c(Q_{\beta})},
\end{cases}
\]
and for all $\beta'\in \{0,1\}^d$ such that ${\rm supp}(\beta) \subseteq {\rm supp}(\beta')$, 
let $\pi_{\beta,\beta'}$ be the projection mapping of $Q_{\beta}$ on $Q_{\beta'}$.
\[
\begin{tikzpicture}
\matrix (m) [matrix of math nodes, row sep=3em, column sep=2.5em, text height=1.5ex, text depth=0.25ex]
{c(Q_{\beta'}) & c(Q_{\beta}) & c(\mathcal{K}_2(|\beta|))\\ }; 
\path[->>,font=\scriptsize] (m-1-1) edge node[auto] {$\check{\pi}_{\beta,\beta'}$} (m-1-2); 
\path[->>,font=\scriptsize] (m-1-2) edge node[auto] {$f_{\beta}$} (m-1-3); 
\end{tikzpicture}
\]
Also, let $\widetilde{K}_{\beta,\beta'}:={\rm Im}(\widetilde{\pi}_{\beta,\beta'}\circ \widetilde{f}_{\beta})$, and 
$g_{\beta,\beta'}: \mathbb{N}^2 \rightarrow \mathbb{N}^{V(Q_{\beta'})}$ given by
\[
g_{\beta,\beta'} (r,t)_{v_{\bf a}}=
\begin{cases}
r & \text{ if } \beta\cdot {\bf a} \text{ is even},\\
t & \text{ if } \beta\cdot {\bf a} \text{ is odd},
\end{cases}
\]
If $\beta'={\bf 1}$, then we simply denote $\pi_{\beta,\beta'}$ by $\pi_{\beta}$, 
$\widetilde{K}_{\beta,\beta'}$ by $\widetilde{K}_{\beta}$, and $g_{\beta,\beta'}$ by $g_{\beta}$.
Note that, if $|\beta|=1$, then $f_{\beta}$ is the identity mapping and if $\beta=\beta'$, 
then $\pi_{\beta,\beta'}$ is the identity mapping.  

\begin{Proposition}\label{main}
Let $d$ be a natural number, $s=V(c(Q_d))\setminus V(Q_d)$, and $\beta,\beta' \in \{0,1\}^d$. 
Then
\begin{description}
\item[$(i)$] $ \widetilde{K}_{\beta}=\{ g_{\beta} (r,t)+ (d-|\beta|){\bf 1}\: | \: 0\leq r,t \leq d \text{ and either } r=|\beta| \text{ or } t=|\beta| \} \triangleleft  SP(c(Q_{d}),s)$,
\item[$(ii)$] $\widetilde{K}_{\beta}$ is generated by $g_{\beta} (d,d-|\beta|)$ with identity $d{\bf 1}$ and $\widetilde{K}_{\beta} \cong {\mathbb Z}_{2|\beta|+1}$,
\item[$(iii)$] $\widetilde{\pi}_{\beta}(SP(c(Q_{\beta}),s))\cap \widetilde{\pi}_{\beta'}(SP(c(Q_{\beta'}),s))=
\widetilde{\pi}_{\beta\odot \beta'}(SP(c(Q_{\beta\odot \beta'}),s))$, 
\end{description}
where $({\bf a} \odot {\bf b})_i={\bf a}_i\cdot {\bf b}_i$ for all $i$.
\end{Proposition}
\begin{proof}
$(i)$ and$(ii)$

By theorems~\ref{uniformhomeo} and ~\ref{cartesian} we get that $\widetilde{\pi}_{\beta}$ and $\widetilde{f}_{\beta}$ are injective homomorphims of groups.
\[
\begin{tikzpicture}
\matrix (m) [matrix of math nodes, row sep=3em, column sep=2.5em, text height=1.5ex, text depth=0.25ex]
{ SP(c(\mathcal{ K}_2(|\beta|)), s) & SP(c(Q_{\beta}),s) & SP(c(Q_d), s) \\};
\path[right hook->,font=\scriptsize]
(m-1-2) edge node[auto] {$\widetilde{\pi}_{\beta}$} 
                        node[auto,swap] {thm.~\ref{cartesian}} (m-1-3);
\path[right hook->,font=\scriptsize]
(m-1-1) edge node[auto,swap] {thm.~\ref{uniformhomeo}}
                        node[auto] {$\widetilde{f}_{\beta}$}  (m-1-2);
\end{tikzpicture}
\]

By theorem~\ref{generators}, $SP(c(Q_{\beta}),s)$ is generated by $\widetilde{f}_{\beta}(|\beta|,0)$ 
with identity $\widetilde{f}_{\beta}(|\beta|,|\beta|)=|\beta|{\bf 1}_{Q_{\beta}}$ for all $\beta\in \{0,1\}^d$.
Moreover, since $\widetilde{K}_{\beta}={\rm Im}(\widetilde{\pi}_{\beta}\circ \widetilde{f}_{\beta})$, then 
\begin{eqnarray*}
\widetilde{K}_{\beta}&=&\{\widetilde{\pi}_{\beta}\circ \widetilde{f}_{\beta}({\bf c}) \,| \, {\bf c} \in SP(c(\mathcal{ K}_2(|\beta|)), s) \}
                                       =\{\widetilde{f}_{\beta}({\bf c})\Box (d-|\beta|){\bf 1}_{Q_{{\bf 1}-\beta}}\, | \, {\bf c} \in SP(c(\mathcal{ K}_2(|\beta|)), s) \}\\
                                    &=&\{ \widetilde{f}_{\beta}((r,t))\Box (d-|\beta|){\bf 1}_{Q_{{\bf 1}-\beta}} \, | \, 0\leq r,t \leq d \text{ and either } r=|\beta| \text{ or } t=|\beta|\}\\
                                    &=&\{ g_{\beta} (r,t)+ (d-|\beta|){\bf 1}\: | \: 0\leq r,t \leq d \text{ and either } r=|\beta| \text{ or } t=|\beta| \} \text{ for all } \beta\in \{0,1\}^d.
\end{eqnarray*}
Thus
\[
\widetilde{K}_{\beta} \cong SP(c(\mathcal{ K}_2(|\beta|)), s)\cong {\mathbb Z}_{2|\beta|+1} \triangleleft SP(c(Q_{d}),s) \text{ for all } \beta\in \{0,1\}^d
\]
is generated by $\widetilde{f}_{\beta}(|\beta|,0)\Box (d-|\beta|){\bf 1}_{Q_{{\bf 1}-\beta}} =
g_{\beta} (d,d-|\beta|)$ and $d{\bf 1}=g_{\bf 1}(d,d)$ is the identity of $\widetilde{K}_{\beta}$.

\begin{figure}[h]\centering
\begin{tikzpicture}[line width=1pt,scale=0.85]
	\tikzstyle{every node}=[inner sep=0pt, minimum width=4pt] 
 	\draw (0,0) {
	+(90:0.82) node[draw, circle, fill=green] (w2) {}
 	+(270:0.82) node[draw, circle, fill=blue] (w3) {}
 	(w2) to [bend right] (w3)
 	(w2) to [bend left] (w3)
 	(w2)+(0,0.35) node {$0$}
 	(w3)+(0,-0.35) node {$2$}
	};
 	\draw (3.7,0) {
	+(45:1.1) node[draw, circle, fill=blue] (u1) {}
 	+(135:1.1) node[draw, circle, fill=green] (u2) {}
 	+(225:1.1) node[draw, circle, fill=blue] (u3) {}
 	+(315:1.1) node[draw, circle, fill=green] (u4) {}
 	(u1) -- (u2) --(u3) -- (u4)-- (u1)
 	(u1)+(0.2,0.2) node {$2$}
 	(u2)+(-0.2,0.2) node {$0$}
 	(u3)+(-0.2,-0.2) node {$2$}
 	(u4)+(0.2,-0.2) node {$0$}
	+(-0.75,0.75) node {\small $Q_{(1,1,0)}$}
	};
	\draw(8,-0.1) {
	+(0,0) node[draw, circle, fill=blue] (v1) {}
	+(0:1.3) node[draw, circle, fill=green] (v2) {}
	+(293:1.015) node[draw, circle, fill=green] (v3) {}
	+(226:1.3) node[draw, circle, fill=blue]  (v4) {}
	(v1)+(0,1.3) node[draw, circle, fill=green] (v5) {}
	(v2)+(0,1.3) node[draw, circle, fill=blue] (v6) {}
	(v3)+(0,1.3) node[draw, circle, fill=blue]  (v7) {}
	(v4)+(0,1.3) node[draw, circle, fill=green] (v8) {}

	[line width=1.5pt] (v1)--(v2)--(v3)--(v4)--(v1)
	[line width=1.5pt] (v5)--(v6)--(v7)--(v8)--(v5)
	[line width=1.5pt] (v1)--(v5) (v2)--(v6) (v3)--(v7) (v4)--(v8)

	(v1)+(-0.3,0.1) node {\small $3$}
	(v2)+(0.3,0) node {\small $1$}
	(v3)+(0.1,-0.3) node {\small $1$}
	(v4)+(-0.15,-0.3) node {\small $3$}
	(v5)+(0,0.35) node {\small $1$}
	(v6)+(0.3,0.3) node {\small $3$}
	(v7)+(0.25,-0.15) node {\small $3$}
	(v8)+(-0.3,0.05) node {\small $1$}
          };
	\draw[right hook->,red,line width=0.5] (0.4,0) to (2.7,0);
	\draw (1.55,-0.3) node[red] {\scriptsize $\widetilde{f}_{(1,1,0)}$};
	\draw[right hook->,line width=0.5] (4.7,0) to (6.7,0);
	\draw (5.85,-0.3) node {\scriptsize $\widetilde{\pi}_{(1,1,0)}$};
	\draw[right hook->, bend left,blue,dashed,line width=0.5] (0.4,0.7) to (7.3,0.7);
	\draw[blue] (3.85,2) node {\scriptsize $\widetilde{\pi}_{(1,1,0)}\circ \widetilde{f}_{(1,1,0)}$};
\end{tikzpicture}
\caption{\small The mappings $\widetilde{f}_{(1,1,0)}$, $\widetilde{\pi}_{(1,1,0)}$, and $\widetilde{\pi}_{(1,1,0)}\circ \widetilde{f}_{(1,1,0)}$ on $Q_3$. }
\end{figure}
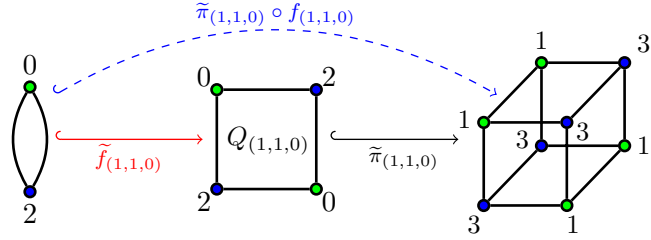	

$(iii)$ Finally, let ${\bf c}\in \widetilde{\pi}_{\beta}(SP(c(Q_{\beta}),s))$.
Since $\widetilde{\pi}_{\beta}({\bf a})={\bf a} \Box {\bf e}$ for all ${\bf a}\in SP(c(Q_{\beta}),s)$, 
where ${\bf e}=(d-|\beta|){\bf 1}$ is the identity of $SP(c(Q_{{\bf 1}-\beta}),s)$.
Then ${\bf c}\in \widetilde{\pi}_{\beta}(SP(c(Q_{\beta}),s))$ if and only if
${\bf c}_{v_{\bf a}}={\bf c}_{v_{\bf b}}$ for all ${\bf a}\odot \beta={\bf b}\odot \beta$.

Now, ${\bf c}\in \widetilde{\pi}_{\beta}(SP(c(Q_{\beta}),s))\cap \widetilde{\pi}_{\beta'}(SP(c(Q_{\beta'}),s))$
if and only if ${\bf c}_{v_{\bf a}}={\bf c}_{v_{\bf b}}$ for all ${\bf a}\odot \beta={\bf b}\odot \beta$, and 
${\bf c}_{v_{\bf a}}={\bf c}_{v_{\bf b}}$ for all ${\bf a}\odot \beta'={\bf b}\odot \beta'$
if and only if ${\bf c}_{v_{\bf a}}={\bf c}_{v_{\bf b}}$ for all ${\bf a}\odot (\beta\odot \beta')={\bf b}\odot (\beta\odot \beta')$
if and only if ${\bf c}\in \widetilde{\pi}_{\beta\odot \beta'}(SP(c(Q_{\beta\odot \beta'}),s))$.
\end{proof}

The next lemma will be useful in order to get the description of the sandpile group of the hypercube. 
\begin{Lemma}{(\cite[Proposition 3.1]{Bai} and \cite[lemma 16]{cartesian})}\label{alpha}
Let $A$ be an abelian group, and let $\alpha$ and $\beta$ be two endomorphisms of $A$
such that $\beta-\alpha=m\cdot 1_{A}$ for some integer $m$.
Then 
\[
Syl_p({\rm coker} \alpha\beta) \cong Syl_p({\rm coker} \alpha \oplus {\rm coker}\beta) 
\]
for all primes $p$ that does not divide $m$.
\end{Lemma}

The next theorem is the main result of this section.

\begin{Theorem}\label{cQd}
Let $k\geq 0$, $d\geq 1$ be natural numbers and let $c_{2k+1}(Q_d)$ be the $2k+1$-cone of the hypercube $Q_d$.
If $s=V(c_{2k+1}(Q_d))\setminus V(Q_d)$, then 
\[
SP(c_{2k+1}(Q_d),s)\cong \bigoplus_{i=0}^{d} \mathbb{Z}_{2i+2k+1}^{\binom{d}{i}}.
\]
Furthermore, $SP(c(Q_{d}),s)=\bigoplus_{\beta \in \{0,1\}^d} \widetilde{K}_{\beta}$.
\end{Theorem}
\begin{proof}
Let $s=V(c_{i}(Q_j))\setminus V(Q_j)$. 
Using elementary row and columns operations invertible over $\mathbb{Z}$ we get an equivalent matrix of $L(c_{2k+1}(Q_{d+1}),s)$.
\begin{eqnarray*}
L(c_{2k+1}(Q_{d+1}),s)
&=& 
\left[
	\begin{array}{cc}
	L(c_{2k+2}(Q_d),s) & -I_{2^d}\\
	-I_{2^d} & L(c_{2k+2}(Q_d),s)
	\end{array}
\right]
\sim
\left[
	\begin{array}{cc}
	-I_{2^d} & L(c_{2k+2}(Q_d),s)\\
	L(c_{2k+2}(Q_d),s) & -I_{2^d}
	\end{array}
\right]\\
&\sim&
\left[
	\begin{array}{cc}
	I_{2^d} & -L(c_{2k+2}(Q_d),s)\\
	0 & L(c_{2k+2}(Q_d),s)^2-I_{2^d}				
	\end{array}
\right]
\sim
\left[
	\begin{array}{cc}
	I_{2^d} & 0\\
	0 & L(c_{2k+1}(Q_d),s)\cdot L(c_{2k+3}(Q_d),s)
	\end{array}
\right]
\end{eqnarray*}
Thus
\begin{eqnarray*}
|SP(c_{2k+1}(Q_{d+1}),s)| & = & |L(c_{2k+1}(Q_{d+1}),s)| = |L(c_{2k+1}(Q_{d}),s)|\cdot |L(c_{2k+3}(Q_{d}),s)|\\
& = & \prod_{i=0}^{d} |L(c_{2k+2i+1}(Q_{1}),s)|^{\binom{d}{i}}=\prod_{i=0}^{d} [(2k+2i+3)(2k+2i+1)]^{\binom{d}{i}}\\
& = & \prod_{i=1}^{d+1} (2k+2i+1)^{\binom{d+1}{i}}.
\end{eqnarray*}

Applying lemma~\ref{alpha} to $A=\mathbb{Z}^{2^{d}}$,
$\alpha=L(c_{2k+1}(Q_{d}),s)$ and $\beta=L(c_{2k+3}(Q_{d}),s)$ ($\beta-\alpha=2I_{2^d}$) we get that
\begin{eqnarray*}
Syl_p(SP(c_{2k+1}(Q_{d+1}),s)) &\cong& Syl_p({\rm coker} \alpha\beta) \cong Syl_p({\rm coker} \alpha \oplus {\rm coker}\beta)\\
&\cong& Syl_p(SP(c_{2k+1}(Q_{d}),s)) \oplus Syl_p(SP(c_{2k+3}(Q_{d}),s)).
\end{eqnarray*}
Therefore, using induction on $d$ and the fact that $(2,|SP(c_{2k+1}(Q_{d+1}),s)|)=1$ we get that 
\[
SP(c_{2k+1}(Q_d),s)\cong \bigoplus_{i=0}^{d} \mathbb{Z}_{2i+2k+1}^{\binom{d}{i}}.
\]

On the other hand, we will use induction on $d$ in order to prove that $SP(c(Q_{d}),s)=\bigoplus_{\beta \in \{0,1\}^d} \widetilde{K}_{\beta}$.
For $d=1$, the result follows from theorem~\ref{generators}.
Now, for all $d\geq 1$, ${\bf 0} \neq \beta \in \{0,1\}^d$, and ${\bf a}\in \{0,1\}^{{\rm supp}(\beta)}$, let 
\[
\Delta_{\beta}({\bf a})=d_{c(Q_{\beta})}(v_{\bf a})e_{v_{\bf a}}-\sum_{{\rm supp}({\bf b})\subseteq {\rm supp}(\beta)}^{{\bf a}- {\bf b}=\pm e_i} e_{v_{\bf b}}
\in \mathbb{Z}^{V(Q_{\beta})}
\]
be the toppling operator of the vertex $v_{\bf a}$ of $(c(Q_{\beta}),s)$, 
$I_{\Delta}(\beta)=\langle \{ \Delta_{\beta}({\bf a}) \,| \, {\bf a}\in \{0,1\}^{{\rm supp}(\beta)}\} \rangle$
be the subgroup generated by the images of the toppling operators of $(c(Q_{\beta}),s)$, and
\[
I_{\Gamma}(\beta)=\langle \{ \Gamma_{\beta',\beta}=g_{\beta',\beta}(|\beta|,|\beta|-|\beta'|)\, | \, {\bf 0} \neq \beta' \in \{0,1\}^{{\rm supp}(\beta)}\} \rangle,
\] 
be the subgroup generated by the generators of $\widetilde{K}_{\beta',\beta}$.
If $\beta={\bf 1}$, we simply denote $\Delta_{\beta}({\bf a})$ by $ \Delta ({\bf a})$, $I_{\Delta}(\beta)$ by $I_{\Delta}$, 
$\Gamma_{\beta',\beta}$ by $\Gamma_{\beta'}$, and $I_{\Gamma}(\beta)$ by $I_{\Gamma}$.

Since, $\bigoplus_{\beta \in \{0,1\}^d} \widetilde{K}_{\beta} \triangleleft SP(c(Q_{d}),s)$ (proposition~\ref{main} $(i)$) and
\[
|SP(c(Q_{d}),s))|=\prod_{i=1}^{d} (2i+1)^{\binom{d}{i}}=\prod_{\beta \in \{0,1\}^d} (2|{\beta}|+1)
=\prod_{\beta \in \{0,1\}^d} |\widetilde{K}_{\beta}|=|\bigoplus_{\beta \in \{0,1\}^d} \widetilde{K}_{\beta}|,
\]
then proving that $SP(c(Q_{d}),s)=\bigoplus_{\beta \in \{0,1\}^d} \widetilde{K}_{\beta}$ 
is equivalent to prove that $\mathbb{Z}^{V(Q_d)}=\langle I_{\Delta}\cup I_{\Gamma}\rangle$.
Hence, we will prove this equivalent form of theorem~\ref{cQd}.

\begin{Theorem}\label{cQdalternativo}
Let $d\geq 1$ be a natural number, then
$\mathbb{Z}^{V(Q_d)}=\langle I_{\Delta}\cup I_{\Gamma}\rangle$
\end{Theorem}
We will use induction on $d$. 
For $d=1$, the result is clear because $\Delta(0)=(2,-1)$, $\Delta(1)=(-1,2)$, 
and by theorem~\ref{generators}, $\Gamma_1=(1,0)$ is a generator of $SP(c(Q_1),s)$ .
Let us assume that the result is true for all hypercubes of dimension less than $d-1$. 

The proof is divided in two steps, first we will prove that $(2d+1)e_{v_{\bf a}} \in \langle I_{\Delta},I_{\Gamma}\rangle$ for all $v_{\bf a}\in V(Q_d)$
and after that we will prove that $d2^{d-1}e_{v_{\bf a}}\in \langle I_{\Delta},I_{\Gamma}\rangle$ for all $v_{\bf a}\in V(Q_d)$.

\begin{Claim}\label{2d+1}
If $d\geq 1$, then $(2d+1)e_{v_{\bf a}} \in \langle I_{\Delta},I_{\Gamma}\rangle$ for all $v_{\bf a}\in V(Q_d)$. 
\end{Claim}
\begin{proof}
We will fix some notation that will be useful in the following.
For all ${\bf 0},{\bf 1} \neq\beta \in \{0,1\}^d$, ${\bf a} \in \{0,1\}^{{\rm supp }(\beta)}$, 
and ${\bf b}\in \{0,1\}^{{\rm supp }({\bf 1}-\beta)}$, let $h_{\beta}^{\bf b}: \{0,1\}^{{\rm supp}(\beta)}\rightarrow \{0,1\}^d$ 
be the mapping given by 
\[
h_{\beta}^{\bf b}({\bf a})_i=
\begin{cases}
{\bf a}_i & \text{ if } i \in {\rm supp}(\beta),\\
{\bf b}_i  & \text{ if }  i \notin {\rm supp}(\beta). 
\end{cases}
\]
Now, let $\beta \in \{0,1\}^d$ with $|\beta|=d-1$.
Since $Q_d\cong Q_{\beta}\Box \mathcal{ K}_2$, then 
\begin{itemize}
\item $\Delta_{\beta}({\bf a})\Box {\bf 0}=\Delta (h^0_{\beta}({\bf a}))+\Delta (h^1_{\beta}({\bf a}))\in I_{\Delta}$ for all ${\bf a}\in \{0,1\}^{{\rm supp}(\beta)}$, 
\item ${\bf 1}_{V(Q_d)}=\sum_{{\bf a}\in \{0,1\}^d} \Delta ({\bf a})\in I_{\Delta}$, and 
\item $\Gamma_{\beta',\beta}\Box {\bf 0}=\Gamma_{\beta'}-{\bf 1}_{V(Q_d)}\in \langle I_{\Delta},I_{\Gamma}\rangle$ 
for all $\beta'$ such that ${\rm supp}(\beta')\subseteq {\rm supp}(\beta)$.
\end{itemize}
Thus, if 
\[
g=\sum_{{\bf a}\in \{0,1\}^{{\rm supp}(\beta)}} z_{\bf a}\Delta_{\beta}({\bf a}) +
\sum_{\beta'\in \{0,1\}^{{\rm supp}(\beta)}} w_{\beta'}\Gamma_{\beta',\beta} \in \langle I_{\Delta}(\beta),I_{\Gamma}(\beta)\rangle,
\]
then $g\Box {\bf 0}\in \langle I_{\Delta},I_{\Gamma}\rangle$.
In particular, by induction hypothesis, $e_{v_{\bf a}} \in \langle I_{\Delta}(\beta),I_{\Gamma}(\beta)\rangle$ for all ${\bf a}\in \{0,1\}^{{\rm supp}(\beta)}$
and therefore 
\[
e_{v_{\bf a}}\Box {\bf 0}=e_{v_{h^0_{\beta}({\bf a})}}+e_{v_{h^1_{\beta}({\bf a})}} \in \langle I_{\Delta},I_{\Gamma}\rangle \text{ for all }{\bf a}\in \{0,1\}^{{\rm supp}(\beta)},
\] 
that is, if $e=v_{\bf a}v_{{\bf a}'}$ is an edge of $Q_d$ and $\chi_e$ is its characteristic vector, then $\chi_e\in \langle I_{\Delta},I_{\Gamma}\rangle$.
Moreover, if $v_{\bf a}$ and $v_{{\bf a}'}$ are vertices of $Q_d$, then 
\[
e_{v_{\bf a}}+(-1)^{{\rm dist}(v_{\bf a},v_{{\bf a}'})} e_{v_{{\bf a}'}} \in \langle I_{\Delta},I_{\Gamma}\rangle,
\]
where ${\rm dist}(v_{\bf a},v_{{\bf a}'})$ is the distance between $v_{{\bf a}}$ and $v_{{\bf a}'}$ in $Q_d$.
Finally, since $\chi_e\in \langle I_{\Delta},I_{\Gamma}\rangle$ for all $e\in E(Q_d)$, then
\[
(2d+1)e_{v_{\bf a}}=\Delta({\bf a})+\sum_{{\bf a}'\in \{0,1\}^d}^{{\bf a}-{\bf a}'=\pm e_i } \chi_{v_{\bf a}v_{\bf a}'} \in \langle 
I_{\Delta},I_{\Gamma}\rangle \text{ for all }v_{\bf a}\in V(Q_d).
\vspace{-10mm}
\] 
\end{proof}

\medskip
Now, we will prove that $d2^{d-1}e_{v_{\bf a}}$ is in $\langle I_{\Delta},I_{\Gamma}\rangle$ for all $v_{\bf a}\in V(Q_d)$.
\begin{Claim}\label{d2d}
If $d\geq 2$, then $d2^{d-1}e_{v_{\bf a}}\in \langle I_{\Delta},I_{\Gamma}\rangle$ for all $v_{\bf a}\in V(Q_d)$.
\end{Claim}
\begin{proof}
Again, we need to fix some notation before beginning with the proof.
For all $n\leq 1$, let $h_{1}: \{0,1\}^{n}\rightarrow \{0,1\}^{n+1}$ be given by
\[
h_{1}({\bf a})_i=
\begin{cases}
{\bf a}_i+1 \, (\text{mod } 2) & \text{ if } i=1,\\
{\bf a}_i & \text{ if } 2\leq i \leq n,\\
1 & \text{ if }  i=n+1, 
\end{cases}
\]
$h_{0}: \{0,1\}^{n}\rightarrow \{0,1\}^{n+1}$ given by $h_{0}({\bf a})=h_{{\bf 1}_n}^0({\bf a})$. 
For $k=1,2$, let $H_k: \mathbb{Z}^{V(Q_n)}\rightarrow \mathbb{Z}^{V(Q_{n+1})}$ be the mapping 
given by
\[
H_k \left( \sum_{v_{\bf a}\in V(Q_n)} z_{v_{\bf a}} e_{v_{\bf a}} \right) =\sum_{v_{\bf a}\in V(Q_n)} z_{v_{\bf a}}e_{v_{h_k({\bf a})}}. 
\] 
Also, for all $d\geq 2$, let 
\[
R_d=H_0\left( \frac{d}{d-1}R_{d-1}\right)+H_1\left( \frac{d}{d-1}R_{d-1}\right)+d2^{d-2}(\chi_{v_{{\bf 0}}v_{e_1}}-\chi_{v_{e_1}v_{h_1({\bf 0})}}),
\]
with $R_2=2(\chi_{v_{(0,0)}v_{(1,0)}} -\chi_{v_{(1,0)}v_{(1,1)}} )$.

\medskip

Actually, we will prove by induction on $d$ that 
$R_d\in \langle I_{\Delta},I_{\Gamma}\rangle$ and $d2^{d-1}e_{v_{\bf 0}}=\Gamma_{\bf 1} +R_d \in \langle I_{\Delta},I_{\Gamma}\rangle$ for all $d\geq 2$.

For $d=2$, clearly $R_2\in \langle I_{\Delta},I_{\Gamma}\rangle$ because $\chi_{e}\in \langle I_{\Delta},I_{\Gamma}\rangle$
for all $e\in E(Q_2)$ and $4e_{v_{(0,0)}}=\Gamma_{(1,1)}+R_2\in \langle I_{\Delta},I_{\Gamma}\rangle$.
Moreover, $4e_{v_{\bf a}}\in \langle I_{\Delta},I_{\Gamma}\rangle$ for all $v_{\bf a}\in V(Q_2)$ because $\chi_{e}$ for all $e\in E(Q_2)$ and $Q_2$ is connected.
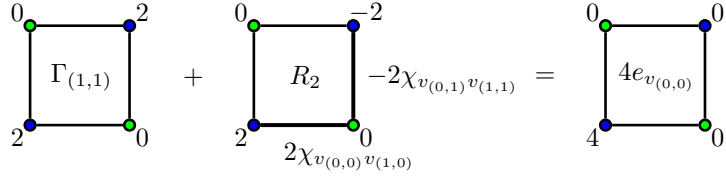
\begin{figure}[h]\centering
\begin{tikzpicture}[line width=1pt,scale=0.85]
	\tikzstyle{every node}=[inner sep=0pt, minimum width=4pt] 
	 \draw (0,0) {
	+(45:1.1) node[draw, circle, fill=blue] (u1) {}
 	+(135:1.1) node[draw, circle, fill=green] (u2) {}
 	+(225:1.1) node[draw, circle, fill=blue] (u3) {}
 	+(315:1.1) node[draw, circle, fill=green] (u4) {}
 	(u1) -- (u2) --(u3) -- (u4)-- (u1)
 	(u1)+(0.2,0.2) node {\small $2$}
 	(u2)+(-0.2,0.2) node {\small $0$}
 	(u3)+(-0.2,-0.2) node {\small $2$}
 	(u4)+(0.2,-0.2) node {\small $0$}
	+(-0.75,0.75) node {\small $\Gamma_{(1,1)}$}
	};
 	\draw (3.5,0) {
	+(45:1.1) node[draw, circle, fill=blue] (u1) {}
 	+(135:1.1) node[draw, circle, fill=green] (u2) {}
 	+(225:1.1) node[draw, circle, fill=blue] (u3) {}
 	+(315:1.1) node[draw, circle, fill=green] (u4) {}
 	(u1) -- (u2) --(u3) -- (u4)-- (u1)
 	(u1)+(0.2,0.2) node {\small $-2$}
 	(u2)+(-0.2,0.2) node {\small $0$}
 	(u3)+(-0.2,-0.2) node {\small $2$}
 	(u4)+(0.2,-0.2) node {\small $0$}
	};
	\draw[line width=1.5pt] (3.5,0) {
	(u3) to (u4) (u4) to  (u1)
	 (u3)+(1.5,-0.5) node {\small $2\chi_{v_{(0,0)}v_{(1,0)}}$ }
 	(u4)+(1.4,0.7) node {\small $-2\chi_{v_{(0,1)}v_{(1,1)}}$ }
	+(-0.75,0.75) node {\small $R_2$}
	+(-2.5,0.75) node {\small $+$}
	+(3,0.75) node {\small $=$}
	};
	\draw (9,0) {
	+(45:1.1) node[draw, circle, fill=blue] (u1) {}
 	+(135:1.1) node[draw, circle, fill=green] (u2) {}
 	+(225:1.1) node[draw, circle, fill=blue] (u3) {}
 	+(315:1.1) node[draw, circle, fill=green] (u4) {}
 	(u1) -- (u2) --(u3) -- (u4)-- (u1)
 	(u1)+(0.2,0.2) node {\small $0$}
 	(u2)+(-0.2,0.2) node {\small $0$}
 	(u3)+(-0.2,-0.2) node {\small $4$}
 	(u4)+(0.2,-0.2) node {\small $0$}
	+(-0.75,0.75) node {\small $4e_{v_{(0,0)}}$}
	};
\end{tikzpicture}
\caption{\small $\Gamma_{(1,1)}+R_2=4e_{v_{(0,0)}}$}
\end{figure}

Let us assume that the result is true for all the natural numbers less or equal than $d-1$.
Since, $v_{h_k({\bf a})}v_{h_k({\bf b})}\in E(Q_{n+1})$ for all $v_{\bf a}v_{\bf b}\in E(Q_n)$, $n\geq 1$, and $k=1,2$,
then $H_k(\frac{d}{d-1}R_{d-1})\in \langle I_{\Delta},I_{\Gamma}\rangle$ for $k=1,2$.
Thus, $R_d=H_0(\frac{d}{d-1}R_{d-1})+H_1(\frac{d}{d-1}R_{d-1})+d2^{d-2}(\chi_{v_{{\bf 0}}v_{e_1}}-\chi_{v_{e_1}v_{h_1({\bf 0})}}) \in \langle I_{\Delta},I_{\Gamma}\rangle$. 
Moreover, since $\Gamma_{\bf 1}=H_0(\frac{d}{d-1}\Gamma_{{\bf 1}_{d-1}})+H_1(\frac{d}{d-1}\Gamma_{{\bf 1}_{d-1}})$, then 
{\small
\begin{eqnarray*}
\Gamma_{\bf 1}+H_0\left(\frac{d}{d-1}R_{d-1}\right)+H_1\left(\frac{d}{d-1}R_{d-1}\right) & = & H_0\left(\frac{d}{d-1}\Gamma_{{\bf 1}_{d-1}}+R_{d-1}\right)+
H_1\left(\frac{d}{d-1}\Gamma_{{\bf 1}_{d-1}}+R_{d-1}\right)\\
& = &  H_0(d2^{d-2}e_{v_{\bf 0}})+H_0(d2^{d-2}e_{v_{\bf 0}})=d2^{d-2} (e_{v_{h_0({\bf 0})}}+e_{v_{h_1({\bf 0})}}),
\end{eqnarray*}
}
and therefore
\[
(d2^{d-1})e_{v_{\bf 0}}=\Gamma_{\bf 1}+R_d\in \langle I_{\Delta},I_{\Gamma}\rangle.
\]
Since $\chi_{e}$ for all $e\in E(Q_d)$ and $Q_d$ is connected, then $d2^{d-1}e_{v_{\bf a}}\in \langle I_{\Delta},I_{\Gamma}\rangle$ for all $v_{\bf a}\in V(Q_d)$.
\end{proof}
Finally, by claims~\ref{2d+1} and ~\ref{d2d}, $(2d+1)e_{v_{\bf a}}\in \langle I_{\Delta},I_{\Gamma}\rangle$ 
and $(d2^{d-1})e_{v_{\bf a}}\in \langle I_{\Delta},I_{\Gamma}\rangle$ for all ${v_{\bf a}}\in V(Q_d)$
and therefore 
\[
e_{v_{\bf a}}\in \langle I_{\Delta},I_{\Gamma}\rangle \text{ for all } {v_{\bf a}}\in V(Q_d)
\] 
because $(2d+1,d2^{d-1})=(2d+1,d)=1$.
\end{proof}

\begin{Remark}
If $n$ is odd, then theorem~\ref{cQd} says that 
\[
SP(c_{n}(Q_d),s)\cong \bigoplus_{i=0}^{d} \mathbb{Z}_{2i+n}^{\binom{d}{i}}.
\]

However, when $n$ is even, this formula is not valid.
For instance, if $n=2$ and $d=2$, then
\[
SP(c_{2}(Q_2),s)\cong \mathbb{Z}_{8}^2\oplus \mathbb{Z}_{3}\not\cong 
\mathbb{Z}_{2}\oplus \mathbb{Z}_{4}^2\oplus \mathbb{Z}_{6}=\bigoplus_{i=0}^{d} \mathbb{Z}_{2i+n}^{\binom{d}{i}}.
\]
\end{Remark}

\begin{Remark}
In~\cite{lorenzini08} is established a close relation between the sandpile group of a graph 
$G$ and the eigenvalues and eigenvectors of their Laplacian matrix.
For instance, Lorenzini~\cite{lorenzini08} proved that if $\lambda$ is an integral eigenvalue of the Laplacian matrix of a graph $G$
and $\mu(\lambda)$ is the maximum number of linear independent eigenvectors associated to $\lambda$, 
then the sandpile group of $G$ contains a subgroup isomorphic to $\mathbb{Z}_{\lambda}^{\mu(\lambda)-1}$, 
see~\cite[Proposition 2.3]{lorenzini08}.

When $G$ is the $n$-cone of the hypercube of dimension $d$ we can use induction on $n$ and $d$ in order to get the eigenvalues of their Laplacian matrix.
More precisely: If $n, d$ are natural numbers and $0\leq i \leq d$, then $\lambda_i=2i+n$ is an eigenvalue of the 
Laplacian matrix of $c_n(Q_d)$ with multiplicity $\binom{d}{i}$.
Using the results in~\cite{lorenzini08} that relate the eigenvalues and eigenvectors of the Laplacian matrix of 
a graph you can get only a partial description of the sandpile group of $c_{2k+1}(Q_d)$.
For instance, using results in~\cite{lorenzini08} you can guarantee only that $\mathbb{Z}_{3}^2 \triangleleft SP(c(Q_2),s) \cong \mathbb{Z}_{3}^2\oplus \mathbb{Z}_{5}$.
The $n$-cone of the hypercube is an example that in general the eigenvalues and eigenvectors of the Laplacian matrix 
are not enough to determine the group structure of the sandpile group.
\end{Remark}

\begin{Remark}
If ${\bf c}_i \in \widetilde{K}_{{\beta}_i}$ for $i=1,2$ with $\beta_1 \odot \beta_2=0$, then ${\bf c}_1\oplus {\bf c}_2={\bf c}_1+{\bf c}_2-d{\bf 1}$.
Also, by theorem~\ref{generators} the generator $g_{\beta}(d, d-|\beta|)$ of $\widetilde{K}_{\beta}$ satisfies that
\[
k\cdot g_{\beta}(d, d-|\beta|)=
\begin{cases}
(d-j,d) & \text{ if } k=2j\leq 2|\beta|,\\
(d,d+j-|\beta|) & \text{ if } k=2j+1 \leq 2|\beta|+1.
\end{cases}
\]
That is, in some cases it is easy to compute the sum of two elements of the sandpile group of $c(Q_d)$.
\end{Remark}

\begin{Example}\label{cq1cq2}
If $d=1$, we have that $SP(c(Q_1))=\mathbb{Z}_3$, $SP(c(Q_1), s)=\{(1,0),(0,1),(1,1)\}$, 
$(1,0)$ and $(0,1)$ are generators $SP(c(Q_1))$, and $(1,1)$ is the identity of $SP(c(Q_1))$.

If $d=2$, we have that $SP(c(Q_2))=\mathbb{Z}_3^2\oplus \mathbb{Z}_5$, $SP(c(Q_2), s)$
is generated by the recurrent configurations $\{\Gamma_{(1,0)}=(2,1,2,1),\Gamma_{(0,1)}=(2,2,1,1),\Gamma_{(1,1)}=(2,0,2,0)\}$, and $\Gamma_{(0,0)}=(2,2,2,2)$ 
is the recurrent configuration that plays the role of the identity in $SP(c(Q_2))$.

\begin{figure}[h]\centering
\begin{tikzpicture}[line width=1.1pt, scale=0.9]
	\tikzstyle{every node}=[inner sep=0pt, minimum width=4pt] 
 	\draw (4,0) {
	+(45:1) node[draw, circle, fill=blue] (v1) {}
 	+(135:1) node[draw, circle, fill=blue] (v2) {}
 	+(225:1) node[draw, circle, fill=blue] (v3) {}
 	+(315:1) node[draw, circle, fill=blue] (v4) {}
 	(v1) -- (v2) --(v3) -- (v4)-- (v1)
          (v1)+(0.2,0.2) node {\small $2$} (v2)+(-0.2,0.2) node {\small $2$} (v3)+(-0.2,-0.2) node {\small $2$} (v4)+(0.2,-0.2) node {\small $2$}
          (v1)+(-0.7,-0.7) node {\small $\Gamma_{(0,0)}$}
          };
          \draw (8,0) {
	+(45:1) node[draw, circle, fill=green] (v1) {}
 	+(135:1) node[draw, circle, fill=blue] (v2) {}
 	+(225:1) node[draw, circle, fill=blue] (v3) {}
 	+(315:1) node[draw, circle, fill=green] (v4) {}
 	(v1) -- (v2) --(v3) -- (v4)-- (v1)
          (v1)+(0.2,0.2) node {\small $1$} (v2)+(-0.2,0.2) node {\small $2$} (v3)+(-0.2,-0.2) node {\small $2$} (v4)+(0.2,-0.2) node {\small $1$}
          (v1)+(-0.7,-0.7) node {\small $\Gamma_{(1,0)}$}
          }; 
 	\draw (12,0) {
	+(45:1) node[draw, circle, fill=green] (v1) {}
 	+(135:1) node[draw, circle, fill=green] (v2) {}
 	+(225:1) node[draw, circle, fill=blue] (v3) {}
 	+(315:1) node[draw, circle, fill=blue] (v4) {}
 	(v1) -- (v2) --(v3) -- (v4)-- (v1)
          (v1)+(0.2,0.2) node {\small $1$} (v2)+(-0.2,0.2) node {\small $1$} (v3)+(-0.2,-0.2) node {\small $2$} (v4)+(0.2,-0.2) node {\small $2$}
          (v1)+(-0.7,-0.7) node {\small $\Gamma_{(0,1)}$}
          };
 	\draw (16,0) {
	+(45:1) node[draw, circle, fill=blue] (v1) {}
 	+(135:1) node[draw, circle, fill=gray] (v2) {}
 	+(225:1) node[draw, circle, fill=blue] (v3) {}
 	+(315:1) node[draw, circle, fill=gray] (v4) {}
 	(v1) -- (v2) --(v3) -- (v4)-- (v1)
          (v1)+(0.2,0.2) node {\small $2$} (v2)+(-0.2,0.2) node {\small $0$} (v3)+(-0.2,-0.2) node {\small $2$} (v4)+(0.2,-0.2) node {\small $0$}
          (v1)+(-0.7,-0.7) node {\small $\Gamma_{(1,1)}$}
          };         
\end{tikzpicture}
\vspace{-2mm}
\caption{\small The identity and the generators of $SP(c(Q_2),s)$}
\end{figure}
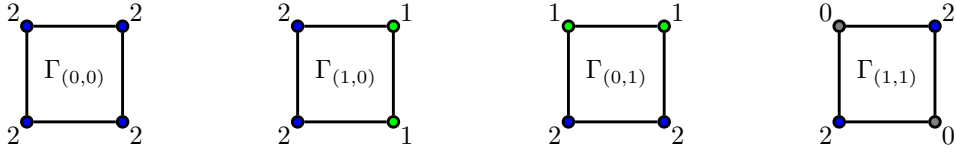
Furthermore, 
\[
\widetilde{K}_{(1,0)}=\{(2,1,2,1),(1,2,1,2),(2,2,2,2)\}
\text{ and } 
\widetilde{K}_{(0,1)}=\{(2,2,1,1),(1,1,2,2),(2,2,2,2)\}
\] 
form two subgroups of $SP(c(Q_2))$ isomorphics to $\mathbb{Z}_3$, and
\[
\widetilde{K}_{(1,1)}=\{(2,0,0,2),(1,2,2,1),(2,1,1,2),(0,2,2,0),(2,2,2,2)\}
\] 
forms one subgroup isomorphic to $\mathbb{Z}_5$ and $SP(c(Q_2),s)=\widetilde{K}_{(1,0)}\oplus \widetilde{K}_{(0,1)}\oplus \widetilde{K}_{(1,1)}$.
\end{Example}

\begin{Remark}
It is clear that 
\[
SP(c_n(Q_0),s)=\{(i)\, | \, 0\leq i \leq n-1\}\cong \mathbb{Z}_{n}
\] 
and $(i)\oplus (j)=(i+j ({\rm mod}\, \, n))$.
Also, it is not difficult to see that ${\bf e}_{(c_n(Q_d),s)}=k_{max} \cdot n{\bf 1}_{2^d}$,
where $k_{max}={\rm max}\{i \, | \, n\cdot i \leq n+d-1\}$ is the identity of $SP(c_n(Q_d),s)$. 
Furthermore, if $n> 1$ and
\[
\hat{K}_{\beta}= \hat{\pi}_{\bf 1}\circ \widetilde{f}_{\beta} (SP(c_n({\mathcal K}_2(|\beta|)),s)) / \hat{\pi}_{\bf 1}(SP(c_n(Q_0),s)),
\]
then $SP(c_n(Q_d),s)= \bigoplus_{\beta \in \{0,1\}^d} \hat{K}_{\beta}$.

For instance, if $d=2$ and $n=3$, then $SP(c_3(Q_2))=\mathbb{Z}_3\oplus \mathbb{Z}_5^2\oplus \mathbb{Z}_7$, $SP(c_3(Q_2), s)$
is generated by the recurrent configurations $\{(2,2,2,2),(3,0,3,0),(0,0,3,3),(0,3,3,0)\}$, and $(3,3,3,3)$ 
is the recurrent configuration that plays the role of the identity in $SP(c_3(Q_2))$.
Furthermore,
\[
\hat{K}_{(0,0)}=\{(2,2,2,2),(4,4,4,4),(3,3,3,3)\}
\]
forms a subgroup of $SP(c(Q_2),s)$ of order $3$,
\begin{eqnarray*}
\hat{K}_{(1,0)}&=&\{(3,0,3,0),(2,1,2,1),(1,2,1,2),(0,3,0,3),(3,3,3,3)\}, \text{ and }\\ 
\hat{K}_{(0,1)}&=&\{(0,0,3,3),(2,2,1,1),(1,1,2,2),(3,3,0,0),(3,3,3,3)\} \nonumber
\end{eqnarray*}
form two subgroups of $SP(c(Q_2),s)$ of order $5$, and 
\[
\hat{K}_{(1,0)}=\{(0,3,3,0),(2,1,1,2),(2,4,4,2),(4,2,2,4),(1,2,2,1),(3,0,0,3),(3,3,3,3)\}
\]
forms a subgroup of $SP(c(Q_2),s)$ of order $7$.

Note that the construction of the $\hat{K}_{\beta}$ subgroups is not canonical, because 
$\hat{\pi}_{\bf 1}((0))=(0,0,0,0)+{\bf e}_{(c_n(Q_d),s)}$ and $\hat{\pi}_{\bf 1}((2))=(2,2,2,2)$.
\end{Remark}

\begin{Remark}
If $IF(d)$ is the number of invariant factors of $SP(c(Q_d))$, then
\[
IF(d)={\rm max}  \left\{ \sum_{\substack{p | 2i+1\\3\leq p\leq 2d+1}} \binom{d}{i} \, | \,  p \text{ is a prime number}\right\}=
\begin{cases} 
6 & \text{ if } d= 4,\\
\sum_{i=0}^{\lfloor \frac{d-1}{3} \rfloor} \binom{d}{1+3i} & \text{ if }  d\neq 4.
\end{cases}
\]
Furthermore, it is not difficult to see that
\[
\lim_{d\rightarrow \infty} \frac{IF(d)}{2^d}=\frac{1}{3}.
\]
\end{Remark}

\noindent {\bf Acknowledgments}

The authors would like to thank V. Reiner, D. Lorenzini and an anonymous referee for their helpful comments.
 
%==================================================================% 
%===========================Bibliografia===========================% 
%==================================================================% 

\end{document}